\documentstyle[12pt]{article}
\begin{document}
\title{Spectral integration and spectral theory
for non-Archimedean Banach spaces.}
\author {S. Ludkovsky, B. Diarra.}
\date{07 August 2001}
\maketitle
\par Permament addresses: 
\par S. Ludkovsky,
Theoretical Department, Institute of General Physics,
\par Str. Vavilov 38, Moscow, 119991 GSP-1, Russia.
\par B. Diarra, Laboratoire de Math\'ematiques Pures, 
\par Complexe Scientifique des C\'ezeaux,
\par 63177, Aubi\`ere, France.
\begin{abstract}
Banach algebras over arbitrary complete non-Archimedean fields
are considered such that operators may be non-analytic.
There are considered different types of Banach spaces over non-Archimedean 
fields. We have determined the spectrum of 
some closed commutative subalgebras of the Banach algebra
${\cal L}(E)$ of the continuous linear operators on a free Banach space
$E$ generated by projectors.
The spectral integration of non-Archimedean Banach algebras is investigated.
For this a spectral measure is defined.
Its several properties are proved. 
The non-Archimedean analog of Stone theorem is proved.
It contains also the case of $C$-algebras $C_{\infty }(X,{\bf K})$.
A particular case of a representation of a $C$-algebra with the help of 
a $L(\hat A,\mu ,{\bf K})$-projection-valued measure is proved.
Spectral theorems for operators and families of commuting linear continuous 
operators on the non-Archimedean Banach space are considered. 
\end{abstract}
\section{Introduction.}
This paper is devoted to the non-Archimedean theory of
spectral integration with the help of the projection-valued measure.
Spectral integration plays very important role in the theory
of Banach algebras, theory of operators and has applications
to the representation theory of groups and algebras in the classical case
of the field of complex numbers $\bf C$ \cite{dun,fell,hew,nai}.
There are also several works about non-Archimedean Banach algebra theory,
which show that there are substantial differences between the non-Archimedean 
and classical cases 
\cite{ber,diar1,diar2,esc1,esc2,esc3,grus,put,put2,roo,vish}.
In the papers \cite{ber,vish} analytic operators over $\bf C_p$ were
considered and the Shnirelman integration of analytic functions
was used, which differs strongly from the non-Archimedean integration
theory related with the measure theory \cite{roo}.
In the non-Archimedean case spectral theory differs from the classical
results of Gelfand-Mazur, because quotients of commutative 
Banach algebras over a field $\bf K$ by maximal ideals 
may be fields $\bf F$, which contain $\bf K$ as a proper subfield
\cite{roo}. In general for each non-Archimedean field $\bf K$
there exists its extension $\bf F$ such that a field ${\bf F}\ne \bf K$
\cite{diar3,sch1}. 
\par Ideals and maximal ideals of non-Archimedean commutative 
$E$-algebras (see \S 5.1.1) 
and $C$-algebras were investigated in \cite{roo,schroo}.
In the works \cite{diar1,diar2} it was shown the failure of the spectral
theory in the non-Archimedean analog of the Hilbert space and
it was shown that even symmetry properties of matrices lead to the enlargment
of the initial field while a diagonalisation procedure.
In the papers \cite{esc1,esc2,esc3} were analysed formulas of spectral
radius and different notions of spectrum and analysed some aspects of
structures of non-Archimedean Banach algebras.
In the book \cite{roo} and references therein general theory of non-Archimedean
Banach algebras and their isomorphisms was considered.
 It was introduced the notion of $C$-algebras in the 
non-Archimedean case apart from the classical
$C^*$-algebras. There are principal differences in the orthogonality
in the Hilbert space over $\bf C$ and orthogonality in the 
non-Archimedean Banach space. Therefore, 
symmetry properties of operators do not play
the same role in the non-Archimedean case as in the classical case.
\par This paper treats another aspects of 
the non-Archimedean algebra theory and theory of operators.
Banach algebras over arbitrary complete non-Archimedean fields
are considered such that operators may be non-analytic.
There are considered different types of Banach spaces over non-Archimedean 
fields. In \S \S 2-4 are considered specific spaces. In \S 5 are 
considered general cases.
\par Let $\bf K$ be a field. 
A non-Archimedean valuation on $\bf K$ is a function $|*|: {\bf K}\to \bf R$
such that:
\par $(1)$ $|x|\ge 0$ for each $x\in \bf K$;
\par $(2)$ $|x|=0$ if and only if $x=0$;
\par $(3)$ $|x+y|\le \max (|x|,|y|)$ for each $x$ and $y\in \bf K$;
\par $(4)$ $|xy|=|x| |y|$ for each $x$ and $y\in \bf K$. \\
The field $\bf K$ is called topologically complete
if it is complete relative to the following metric:
$\rho (x,y)=|x-y|$  for each $x$ and $y\in \bf K$.
A topological vector space $E$ over $\bf  K$ with  the non-Archimedean 
valuation may have a norm $\| * \| $ such that its restriction on each
one dimensional subspace over $\bf K$ coincides with the valuation
$|*|$. If $E$ is complete relative to such norm $\| * \| $,
then it is called the Banach space. 
Such fields and topological vector spaces are called non-Archimedean.
An algebra $X$ over $\bf K$ is called Banach, if it is a Banach space
as a topological vector space and the multipilication in it
is continuous such that $\| xy \| \le \| x \| \| y \| $
for each $x$ and $y$ in $X$.
A finite or infinite sequence $(x_j: j\in \Lambda )$ 
of elements in a normed space $E$ is called orthogonal,
if $\| \sum_{j\in \Lambda }\alpha _jx_j \|
=\max ( \| \alpha _j x_j \| : j\in \Lambda )$ 
for each $\alpha _j\in \bf K$ for which $\lim_j\alpha _jx_j=0$.
We consider the infinite
topologically complete field $\bf K$ with the nontrivial 
non-Archimedean valuation. 
\par A non-Archimedean Banach space
$E$ is said to be free
if there exists a family $(e_j: {j \in I}) \subset E$ such that any
element $x \in E$ can be
written in the form of convergent sum $x=\displaystyle\sum_{j \in
I}x_je_j$, i.e.,
$\displaystyle\lim_{j \in I}x_je_j =0$,
and $\Vert x \Vert=\displaystyle\sup_{j \in I}\vert x_j \vert \Vert
e_j \Vert$ (see \S 2). In \S 3 ultrametric Hilbert spaces are considered.
In \S 4 we have determined the spectrum of 
some closed commutative subalgebras of the Banach algebra
${\cal L}(E)$ of the continuous linear operators of $E$ generated by
projectors.
\par \S 5 is devoted to the spectral integration.
We introduce another definition of $E$-algebras in \S 5.1 apart from 
\cite{schroo}. In Propositions 5.3 and 5.4 
we have proved that they are contained 
in the class of $E$-algebras and $C$-algebras considered in
\cite{roo,schroo}. In \S 5.2 a spectral measure is defined.
In \S \S 5.6-5.10 its several properties are proved. In \S 5.12 
the non-Archimedean analog of Stone theorem is proved.
It contains also the case of $C$-algebras $C_{\infty }(X,{\bf K})$.
A particular case of a representation of a $C$-algebra with the help of 
$L(\hat A,\mu ,{\bf K})$-projection-valued measure is proved in Theorem 5.15.
Spectral theorems for operators and families of commuting linear continuous 
operators on a non-Archimedean Banach space are considered in 
\S \S 5.17, 5.18.
\section{Free Banach spaces}
{\bf 1.} Let $E$ be the free Banach space with an orthogonal base
$(e_j: {j \in I} )$. The
topological dual $E'$ of $E$ is a Banach space with respect to the norm
defined for
$x' \in E'$ by $\Vert x' \Vert=\displaystyle\sup_{x \not=0}
{\vert <x',x> \vert \over \Vert x \Vert}$.
For $x' \in  E'$ and $y \in E$, one defines an element  $(x'\otimes
y)$ of the Banach algebra of continuous linear operators
${\cal L}(E)$ on the space $E$ by setting for $x \in E$,
$(x'\otimes y)(x)=<x',x>y$ with norm $\Vert x'\otimes y \Vert=\Vert
x' \Vert \Vert y \Vert$.
If $E$ is a free Banach space with base $(e_j: {j \in I})$,
any $u \in {\cal L}(E)$ can be
written as a  pointwise convergent sum $u=\displaystyle\sum_{(i,j) \in
I\times I}\alpha_{ij}
e'_j \otimes e_i$. Hence $ \displaystyle\lim_{i \in I}\alpha_{ij}e_i=
0$ for each $j \in I$.
Moreover $\Vert u \Vert= \displaystyle\sup_{i,j}\vert \alpha_{ij} \vert
\Vert e'_j \Vert \Vert e_i \Vert$. Notice that $\Vert e'_j \Vert={1
\over \Vert e_j \Vert}$.
Let ${\cal L}_0(E)=\{ u: \quad u=\displaystyle\sum_{(i,j) \in I
\times I}\alpha_{ij}
e'_j \otimes e_i \in {\cal L}(E); \quad \displaystyle\lim_{j
\in I}\alpha_{ij}e'_j=
0$ for each $i \in I \}$.
\par {\bf 1.1. Theorem.} {\it ${\cal L}_0(E)$ is a closed
subalgebra in ${\cal L}(E)$ with the unit element of ${\cal L}(E)$}.
\par {\bf Proof.} Let $ u, v \in {\cal L}_0(E), 
u=\displaystyle\sum_{(i,j) \in I \times
I}\alpha_{ij}
e'_j \otimes e_i$  and $v=\displaystyle\sum_{(i,j) \in I \times
I}\beta_{ij} e'_j \otimes e_i$,
then $\displaystyle\lim_{i \in I}\alpha_{ij}e_i=0=\displaystyle\lim_{i
\in I}\beta_{ij}e_i$
for each $j \in I$, and   $\displaystyle\lim_{i \in
I}\alpha_{ij}e_i=0=\displaystyle\lim_{i \in
I}\beta_{ij}e_i$ for each $j \in I$.
One has $u \circ v=\displaystyle\sum_{(i,j) \in I \times I}
(\displaystyle\sum_{k \in I}\alpha_{ik}\beta_{kj})e'_j \otimes e_i$.
Let $i \in I, \quad
\displaystyle\lim_{k \in I}\alpha_{ik}e'_k=0$, that is,
for each $\varepsilon > 0 $, there exists
$J_\varepsilon(i)$ a finite subset of $I$ such that
for each $k \not\in J_{\varepsilon}(i),
\Vert \alpha_{ik}e'_k \Vert < \varepsilon$.
Hence $ \Vert (\displaystyle\sum_{k \in I}\alpha_{ik}\beta_{kj})e'_j
\Vert=\Vert \displaystyle \sum_{k \in J_\varepsilon(i)
}(\alpha_{ik}\beta_{kj})e'_j +$ 
$\displaystyle \sum_{k
\not\in J_\varepsilon (i)}(\alpha _{ik}\beta _{kj})e'_j \Vert
\leq
\max {(\displaystyle \max_{k\in J_\varepsilon (i)} \Vert \alpha _{ik}e'_k
\Vert \Vert e_i \Vert
\Vert \beta_{kj}e_k\Vert \Vert e'_i \Vert \Vert e'_j\Vert ,
\displaystyle \sup_{k\not\in  J_\varepsilon (i)}
\Vert \alpha _{ik}\beta _{kj}e'_j\Vert )}
\leq $ \\
$\max {(\Vert u \Vert \displaystyle  \max_
{k  \in J_\varepsilon (i)} \Vert \beta _{kj}e'_j \Vert \Vert e_k\Vert
\Vert e'_i \Vert ,
\varepsilon \Vert v \Vert \Vert e'_i \Vert )}$.
Since $\displaystyle \lim_{j \in I}\Vert \beta _{kj} e_j \Vert=0$
for each $k \in J_{\varepsilon }(i)$,
one has $\displaystyle \lim_{j \in I}\Vert (\displaystyle \sum_{k \in
I}\alpha _{ik}\beta _{kj})e'_j
\Vert=0$ for each $i \in I$, therefore $u \circ v \in {\cal L}_0(E)$.
The identity  map $id$ being given by $id=\displaystyle \sum_{i \in
I}e'_i \otimes e_i$, one
has $\alpha _{ii}=1$ and $\alpha _{ij}=0$ if $i \not= j$ .Therefore
$\displaystyle \lim_{i}
\alpha _{ij}e_i=0$ for each $j \in I$, and  $\displaystyle \lim_{j}
\alpha _{ij}e'_j=0$ for each $i \in I$. Hence $id \in {\cal L}_0(E)$.
Let $u=\displaystyle \sum_{(i,j) \in I \times I}\alpha _{ij} e'_j
\otimes e_i$ be in the
closure of ${\cal L}_0(E)$. For all $\varepsilon > 0$, there exists
$u_\varepsilon=\displaystyle\sum_{(i,j) \in I \times
I}\alpha_{ij}(\varepsilon) e'_j \otimes e_i \in
{\cal L}_0(E)$ such that $\Vert u-u_\varepsilon
\Vert=\displaystyle\sup_{i,j}\vert\alpha_{ij}-\alpha_{ij}(\varepsilon)\vert
\Vert e'_j \Vert
\Vert e_i \Vert < \varepsilon$. Hence for all $ i, j \in I$, one has
$\vert
\alpha_{ij} \vert
\Vert e'_j \Vert \Vert e_i \Vert  \leq  \max(\varepsilon, \vert
\alpha_{ij}(\varepsilon) \vert \Vert e'_j
\Vert \Vert e_i \Vert )$. One obtains
$\displaystyle\lim_{i}\Vert\alpha_{ij} e_i \Vert=0$
for each $j \in I$  and  $\displaystyle\lim_{j}\Vert\alpha_{ij} e'_j
\Vert=0$ for each $i \in I$.
Therefore  $u \in {\cal L}_0(E)$ and  ${\cal L}_0(E)$ is closed.
\par {\bf 2.} 
Suppose that the orthogonal basis is orthonormal, i.e., $\Vert
e_j \Vert=1$ for each 
$j \in I$. Then $u=\displaystyle\sum_{(i,j) \in I \times I}\alpha_{ij}
e'_j\otimes e_i \in
{\cal L}_0(E)$, if and only if $\displaystyle\lim_{i}\alpha_{ij}=0$
for each $j \in I$  and  $\displaystyle\lim_{j}\alpha_{ij} =0$
for each $i \in I$.
Setting for $u =\displaystyle\sum_{(i,j) \in I \times
I}\alpha_{ij} e'_j \otimes e_i \in  {\cal L}_0(E)$,
$u^*=\displaystyle\sum_{(i,j) \in I \times I}\alpha_{ji} e'_j \otimes
e_i$, one sees that $u^* \in  {\cal L}_0(E)$ , called the adjoint
of $u$. One verifies easily the following.
\par {\bf 2.1. Proposition.} {\it An element $ u \in {\cal L} ( E)$ has
an adjoint $ u^*$ if and only if $ u \in {\cal L}_0 ( E)$ .
Let $ u, v \in {\cal L}_0 ( E) , \lambda \in K.$ One  has $ ( u +
\lambda v )^* = u^* + \lambda v^* ;
( u \circ v ) ^* = v^* \circ u^* ; \quad  u^{**} = u $.  \quad Moreover
$ \Vert u^* \Vert = \Vert u \Vert$.}
\par As usually, one says that $u \in {\cal L}_0 ( E) $ is normal
(respectively unitary)
if $u \circ u^*=u^* \circ u$  (respectively $u \circ u^*= id =u^* \circ u)$.
And $u$ is self-adjoint if $u=u^*$, this is
equivalent here to say that the matrix of $u$ is symmetric.
\par {\bf 3. Note.} $(i)$. One  has  $\Vert u \Vert =\Vert u^* \Vert $.
However , in general
$\Vert u \circ u^*  \Vert \not= \Vert u \Vert^2$.
For example, if $I$ is the set of positive integers, and
$E$ with orthogonal base
$(e_n: {n \geq 1}),$ let  $ a, b \in K$. The operator $u$ defined by
$u(e_1)=ae_1+be_2 , u(e_2)=
be_1-ae_2 , u(e_3)= ce_3$, and $u(e_n)=0$ for $n \geq 4$. One sees that
$u$ is self-adjoint.
If $i=\sqrt{-1} \in K$; then taking $b=ia$ and  $\vert c \vert  <
\vert a \vert $, one sees that
$\Vert u^2 \Vert =\vert c \vert ^2  < \vert a \vert ^2=\Vert u \Vert ^2$.
\par $(ii)$. It should be interesting to characterize the
elements of ${\cal L}_0 ( E) $
that are normal, unitary.
Considering,  whenever the base of $E$ is orthonormal,the bilinear
form  $f$  on $E$  defined
by $f(x,y)=\displaystyle\sum_{ i \in I}x_iy_i$ ,one obtains that the
above definition of an
adjoint $u^*$  of an element $u \in {\cal L}_0 ( E) $ is equivalent to
say that $f(u(x),y)=f(x,u^*(y))$ for each $x$ and $y \in E$. 
In fact, here the adjoint of an operator is its
transposition. This example is related to ultrametric Hilbert spaces.
\section{Ultrametric Hilbert spaces.}
For the so called ultrametric Hilbert spaces one can also define
the adjoint of an operator
with respect to an appropriate bilinear symmetric form.
\par {\bf 1. Remark and Definition.}
H. Ochsenius and W.H. Schikhof write in \cite{ocsc} 
{\it "as a slogan : There
are no $p$-adic Hilbert spaces".}
Nevertheless we shall give a definition of $p$-adic Hilbert spaces [cf.
for example, \cite{ko1,ko2} for some fields with infinite rank valuation].
Let $ \omega = (\omega _i)_{i \geq 0} $ be a sequence of {\it
non-zero} elements of
$\bf K$. Let us consider the free Banach space $E_{\omega }=c_0({\bf N},
{\bf K},( \vert \omega _i \vert ^{1\over2})_{i\ge 0})= 
\{ x: \quad x =(x_i)_{i \ge 0} \subset {\bf K};
\displaystyle \lim_{i \rightarrow +\infty } \vert 
x_i \vert \vert \omega _i \vert ^{1\over 2}=0 \} $.
Then $x=(x_i)_{i \ge 0} \in E_{\omega }\Longleftrightarrow
\displaystyle  \lim_{i \rightarrow
+\infty }x^2_i \omega _i=0$. Setting $e_i=(\delta _{i,j})_{j \ge 0} $
(Kronecker symbol), one has
that $(e_i: i \ge 0)$ is an orthogonal base of $E_{\omega }: \quad
\forall x \in E_{\omega },\quad x=$
$\displaystyle \sum_{i \ge 0} x_i e_i$ and
$\Vert x \Vert = \displaystyle \sup_{i \ge 0} \vert x_i \vert \Vert e_i
\Vert =$ 
$\displaystyle \sup_{i \ge 0} \vert x_i \vert \vert \omega _i
\vert ^{1\over 2}$, in  particular,
$ \Vert e_i \Vert = \vert \omega _i \vert ^{1\over 2}$ for each 
$i \ge 0.$
Let $ f_{\omega }: E_{\omega }\times E_{\omega }\rightarrow \bf K$ be defined by
$f_{\omega }(x,y) =
\displaystyle \sum_{i \ge 0} \omega _i x_i y_i$. It is readily seen
that $f_{\omega }$ is a bilinear
symmetric form on $E_{\omega }$, with $\vert f_{\omega }(x,y) \vert \le
\Vert x \Vert \Vert y \Vert $, 
i.e., the bilinear form $f_{\omega }$ is continuous.
Moreover $f_{\omega }$ is non-degenerate, i.e. $f_{\omega }(x,y)= 0$
for each $y \in  E_{\omega }
\Longrightarrow x = 0$. Furthermore $f_{\omega }(x,x)= \displaystyle
\sum_{i \ge 0} \omega _i x^2_i$ and
$f_{\omega }(e_i,e_j)= \omega _i \delta _{i,j}$ for $i$ and $j \ge 0.$
The space $E_{\omega }$ is called a $p$-adic Hilbert space.
\par {\bf 2. Note.} $(i)$. 
It may happen that $ \vert f_{\omega }(x,x)\vert <
\Vert x \Vert ^2$ for
some $ x \in E_{\omega }$ and even worse, $E_{\omega }$ contains isotropic
elements $ x \not= 0$ , i.e.,  $f_{\omega }(x,x)= 0$
\par $(ii)$. Let $V$ be a subspace of $E_{\omega }$ and 
$V^{\perp }= \{ x \in E_{\omega }: \quad f_{\omega }(x,y)=0, 
\forall y \in V \} $.
The fundamental property on subspaces of the classical Hilbert space 
$H: V =V^{\perp \perp } \Rightarrow V \oplus  V^{\perp }=H$ 
fails to be true in the $p$-adic case. This explains the
claim of H. Ochsenius and  W.H. Schikhof.
\par {\bf 3. Remark.} A free Banach space $E$ with an orthogonal
base $ ( e_i : i
\ge 0)$ can be given a structure of a $p$-adic Hilbert space if and only
if there exists $ (
\omega_i: {i \ge 0}) \subset \bf K$ such that $\Vert e_i \Vert = \vert
\omega_i\vert^{1\over2} $ for each $i \ge 0$.
Furthermore if $\bf K$ contains a square of any of its element,
then any $p$-adic Hilbert is isomorphic, in a natural way,  to the
space $c_0 ({\bf N},{\bf K})$.
\par {\bf 4. Note.}
Let $ u , v \in {\cal L} ( E_\omega ) $; one has  $ u = \displaystyle
\sum_{i , j } \alpha_{ij} e'_j
\otimes e_i $ and $ v = \displaystyle\sum_{ i , j} \beta_{ij} e'_j
\otimes e_i $ with
$\displaystyle \lim_{i \rightarrow + \infty} \vert \alpha _{ij} \vert
\vert \omega _i \vert ^{1\over 2}=0= \displaystyle \lim_{i \rightarrow
+\infty}
\vert \beta _{ij} \vert \vert \omega _i \vert ^{1\over 2}$
for each $j \geq 0$. Furthermore the norm
of $u\in {\cal L} ( E_{\omega })$ is given  by  $\Vert u \Vert
=\displaystyle \sup_{i,j}
{\vert \omega _i \vert ^{1\over 2} \vert \alpha _{ij} \vert\over \vert
\omega _j \vert ^{1\over2}}$.
\par The operator $v$ is said to be an 
adjoint of $u$ with respect to $f_{\omega }\Longleftrightarrow
f_{\omega }(u(x),y)=f_{\omega }(x,v(y)),$ for all $x, y \in E_{\omega }$.
Since $f_{\omega }$ is symmetric $u$ is an adjoint of $v$.   
\par Since $f_{\omega }$ is non degenerate, 
if an operator $u$ has an adjoint,
this adjoint is unique and will be denoted by $ u^*$.
Since $(e_i: i \ge 0)$ is an orthogonal base of $ E_{\omega }$, one
has that $v $ is an
adjoint of $u$ if and only if $ f_{\omega }((u(e_i),e_j) = 
f_{\omega }(e_i,v(e_j))$ for each $i$ and
$j \geq 0$. That is, $ f_{\omega }(\displaystyle \sum_{ k \ge 0}
\alpha _{ki} e_{k},e_j) =\alpha _{ji}\omega _j =$
$f_{\omega }(e_i,\displaystyle \sum_{k\geq 0} \beta _{kj}e_k) =
\beta _{ij}\omega _i, \quad \forall i, j \ge 0
\Longleftrightarrow \beta _{ij}= \omega ^{-1}_i \omega _j \alpha _{ji}, 
\quad \forall i, j \ge 0$.
Furthermore one must have $\displaystyle \lim_{ i \rightarrow +\infty}
\vert \beta _{ij} \vert \vert
\omega _i \vert ^{1\over 2} = 0$ for each $j \ge 0$, that is,
$$\displaystyle \lim_{i \rightarrow
+\infty} \vert \omega _i \vert ^{1\over 2} \vert \omega ^{-1}_j\vert \vert
\omega _j \vert \vert
\alpha _{ji} \vert = \vert \omega _j \vert \displaystyle \lim_{i
\rightarrow +\infty} \vert \omega _i
\vert ^{-{1\over 2}} \vert \alpha _{ji}\vert = 0, \quad \forall j \ge 0.$$
Hence $\displaystyle \lim_{i \rightarrow + \infty }
\vert \omega _i \vert ^{-{1\over 2}} \vert \alpha _{ji} \vert = 0$
for each $j \ge 0$. We have proved the following.
\par {\bf 5. Theorem.} {\it Let $(\omega _i)_{i \ge 0} \subset \bf K^* $ 
and $E_{\omega }=c_0({\bf N},{\bf K}, 
(\vert \omega _i \vert ^{1\over 2})_{i \ge 0})$ be the
$p$-adic Hilbert space associated with $\omega $.
Let $u=\displaystyle \sum_{i,j} \alpha _{ij} {e'}_j \otimes e_i \in
{\cal L} ( E_{\omega })$. Then
$u$ has an adjoint $v = u^* \in {\cal L}( E_{\omega }) $ if and only if
$\lim_{ j \rightarrow
+\infty} \vert \omega _j \vert ^{-{1\over 2}} \vert \alpha _{ij}\vert = 0$
for each $i \ge 0$. In
this condition, $ u^* = \displaystyle \sum_{i,j} \omega ^{-1}_i \omega _j
\alpha _{ji} {e'}_j \otimes e_i$.}
\par It follows from the above theorem that not any continuous linear
operator of $E_{\omega }$ has
an adjoint: it is another difference  with classical Hilbert spaces.
Let ${\cal L}_0 ( E_\omega ) = \{ u: \quad u = 
\displaystyle\sum_{i \ge 0} \sum_{
j \ge 0} \alpha_{ij} e'_j
\otimes e_i \in {\cal L} ( E_\omega ); \quad
\displaystyle\lim_{ j \rightarrow +\infty} \vert \omega_j
\vert^{-{1\over 2}} \vert \alpha_{ij} \vert = 0 , \ \forall i \ge 0\}$.
Let us remind that $ u =
\displaystyle \sum_{ i,j} \alpha_{ij} e'_j \otimes e_i \in {\cal L} (
E_\omega )$ is equivalent to
$\displaystyle\lim_{i \rightarrow +\infty} \vert \omega_i \vert ^{1\over
2} \vert \alpha_{ij} \vert = 0$
for each $j \ge 0$.
It is readily seen, as in Theorem $1$, that ${\cal L}_0 ( E_\omega )$ is
a closed unitary subalgebra of ${\cal L} ( E_\omega )$.
\par {\bf 6. Corollary.} {\it An element $ u \in {\cal L} ( E_\omega)$
has an adjoint $
u^*$ if and only if $ u \in {\cal L}_0 ( E_\omega )$ .
Let $ u, v \in {\cal L}_0 ( E_\omega) , \lambda \in K.$ One  has $ ( u +
\lambda v )^* = u^* + \lambda v^* ;
( u \circ v ) ^* = v^* \circ u^* ; \quad  u^{**} = u $.  \quad Moreover $
\Vert u^* \Vert = \Vert u \Vert$.}
\par {\bf Proof.} We only prove  $\Vert u^* \Vert = \Vert u \Vert $. Since
for $ u = \displaystyle
\sum_{i,j} \alpha_{ij} e'_j \otimes e_i \in {\cal L}_0 ( E_\omega )$ one
has $\Vert u \Vert
=\displaystyle\sup_{i,j} {\vert\omega_i \vert^{1\over 2} \vert
\alpha_{ij} \vert\over \vert \omega_j
\vert^{1\over2}}$ and $ u^* = \displaystyle \sum_{ i,j} \omega_j
\omega_i^{-1} \alpha_{ji} e'_j
\otimes e_i $ , one obtains  $ \Vert u^* \Vert = \displaystyle
\sup_{i,j} {\vert \omega_i 
\vert^{1\over 2}\over \vert \omega_j\vert^{1\over2}} \vert \omega_j
\vert \vert\omega_i^{-1} \vert
\vert \alpha _{ij}\vert =\displaystyle \sup_{ i,j} {\vert \omega _j
\vert ^{1\over2 }\over
\vert \omega_i \vert ^{1\over2 }}\vert \alpha _{ji}\vert = \Vert u \Vert $.
\par {\bf 7. Remark.} $(i)$. $ u = \displaystyle \sum_{i,j} \alpha_{ij}
e'_j \otimes e_i \in
{\cal L}_0  ( E_\omega) $ is self-adjoint, i.e., $ u =    u^* $ if and
only if $ \alpha_{ji} =
\omega_i \omega_j^{-1}\alpha_{ij}$,  for each $i \ge 0$ and each $j \ge 0.$
\par $(ii)$. Examples of self-adjoint operators on
ultrametric Hilbert spaces and
study of their spectrum are given in \cite{ack,abpk,khrb,diar2}.
\section{Closed subalgebras generated by projectors.}
{\bf 1.} Let $J$ be a subset of $I$  and $E$  be a free Banach space
with orthogonal basis
$(e_j: {j \in I})$. The  linear operator $p_J=\displaystyle\sum_{i \in
J}e'_i\otimes e_i$ of $E$ belongs to  ${\cal L}_0(E)$. 
Let ${\cal D}=\{ u: \quad 
u=\displaystyle \sum_{i \in I}\lambda _ie'_i\otimes e_i \in {\cal L}_0(E);
\quad \displaystyle \sup_{i \in I}\vert \lambda_i \vert <
+\infty \} $. It is clear that
${\cal D}$ is isometrically isomorphic to the the algebra of bounded
families $\ell ^{\infty}(I,{\bf K})$.
Let $Hom_{alg}({\cal D}, {\bf K})$ denotes a family of all
algebra homomorphisms of $\cal D$ into $\bf K$.
Consider the spectrum ${\cal X}({\cal D})=Hom_{alg}({\cal D}, {\bf K})$ 
in a topology inherited from the Tihonov topology of the product 
${\bf K}^{\cal  D}$ of copies of $\bf K$.
\par {\bf 1.1. Proposition.} {\it  $(i)$.  
An element $u=\displaystyle \sum_{i \in
I}\lambda _ie'_i\otimes e_i \in {\cal D}$ is
an idempotent if and only if there exists $J \subset I$ such that
$u=p_J$. 
\par $(ii)$. The spectrum 
${\cal X}({\cal D})$ is homeomorphic to the subset of
ultrafilters on $I : \Phi_c=\{ {\cal U}$: $\cal U$ is an
ultrafilter on $I$, such that for all 
\par $u=\displaystyle \sum_{i \in I}\lambda _ie'_i\otimes e_i \in {\cal D}$,
the limit
$\displaystyle \lim_{ {\cal U}}\lambda_i$ exists in $\bf K$ \} .}
\par {\bf Proof.} $(i)$. Let  
$u=\displaystyle \sum_{i \in I}\lambda _ie'_i\otimes e_i $;
then $u \circ u =u$
if and only if $\displaystyle \sum_{i \in I}{\lambda _i}^2e'_i\otimes
e_i=\displaystyle \sum_{i \in
I}\lambda _ie'_i\otimes e_i $, if and only if $\lambda _i^2=\lambda _i$
for each $i \in I,$
if and only if $\lambda _i=0$  or  $\lambda _i=1$ . Setting 
$J=\{ i: \quad i \in I; \quad \lambda _i=1 \} ,$  one has $u=p_J$.  
\par $(ii)$. Let $\chi $ be a character of ${\cal D}$, that is an algebra
homomorphism (necessary
continuous)  of  ${\cal D}$ into $\bf K$.
For all $J , L  \subset I$  one has $p_J \circ p_L=p_{J \cap L}$,
hence  $p_J \circ
p_{J^c}= p_{\emptyset}=0$,where $J^c=I \setminus J$. Furthermore
$\chi (p_J)= \chi (p_J)\chi (p_J)$ implies that  $\chi (p_J)=0$ or $1$.
Let ${\cal U}_{\chi }= \{ J: \quad J \subset I; \quad  \chi (p_J)=1
\}$. This family of subsets
is an ultrafilter. Indeed, $\emptyset \not\in {\cal U}_{\chi }$ . If $J
\subset L$ with $J \in
{\cal U}_{\chi}$,one has $1=\chi (p_J)= \chi (p_{J \cap
L})=\chi (p_J)\chi (p_L)=\chi (p_L)$, hence
$L \in {\cal U}_{\chi}$.
On the other hand, for $J \subset I$,  one has  $1_E=p_J + p_{J^c}$,
and $1=\chi (1_E)=\chi (p_J)+\chi (p_{J^c})$ with
$\chi (p_J) =1$, or $0$ and $\chi (p_{J^c}) =1$, or $0$. If $\chi (p_J)
=1$, one has $\chi (p_{J^c})
=0$, and if $\chi (p_{J^c})=1$, one has $\chi (p_J) =0$. Hence $J \in
{\cal U}_{\chi}$ or $J^c \in {\cal U}_{\chi}$.
Let $u=\displaystyle\sum_{i \in I}\lambda_ie'_i\otimes e_i  \in
{\cal D}$. Put $\chi (u)=
\lambda  \in  K$;  then for all $J \in  {\cal U}_{\chi}$, one has $
\chi (up_J)=\chi (u)=\lambda
\chi (p_J)$. Therefore $\chi (up_J-\lambda p_J)=0$, i.e., $up_J-\lambda
p_J \in  \ker \chi $.
Set $\phi_{{\cal U}_\chi }(u)=\displaystyle\lim_{{\cal
U}_\chi }\vert \lambda_i \vert $. It
is  well known and readily seen that $\phi_{{\cal U}_\chi }$ is a
multiplicative semi-norm on
${\cal D}$ and that $\ker \phi _{{\cal U}_{\chi }}=
\{ u: \quad u \in {\cal D}; \quad 
\phi _{{\cal U}_{\chi }}(u)=0 \} $ is a maximal ideal of ${\cal D}$, since
${\cal D}$ is isomorphic to $\ell ^{\infty }(I,{\bf K})$.
On the other hand $\vert \chi (up_J) \vert \le \Vert up_J \Vert=
\displaystyle \sup_{i \in J}\vert \lambda_i \vert $ for each $J \in {\cal
U}_{\chi }$. It follows
that $\vert \chi (u) \vert=\vert \chi (up_J) \vert \le
\displaystyle \inf_{J \in {\cal U}_\chi }
\displaystyle \sup_{i \in J}\vert \lambda_i \vert=\phi _{{\cal
U}_\chi }(u)$. Hence, one has
$\ker \phi _{{\cal U}_{\chi }} \subset \ker \chi $ and $\ker \phi _{{\cal
U}_\chi } = \ker \chi $.
Let  $J \in {\cal U}_{\chi }$, one deduces from $up_J-\lambda p_J
\in \ker \chi =\ker \phi  _{{\cal U}_{\chi }}$, 
that $0= \phi_{{\cal U}_{\chi }}(up_J-\lambda
p_J)=
\displaystyle \lim_{{\cal U}_{\chi }}\vert \lambda _i-\lambda \vert $. It
follows that $\displaystyle \lim_{{\cal U}_{\chi }}\lambda _i= \lambda $ 
exists in $\bf K$.
Moreover, $\chi (u)=\lambda = \displaystyle \lim_{{\cal U}_{\chi }}\vert
\lambda _i \vert $, and one sees that $\chi =\chi _{{\cal U}_{\chi }}$.
Reciprocally, if ${\cal U}$ is an ultrafilter on $I$ such
that for all
$u=\displaystyle \sum_{i \in I}\lambda _ie'_i\otimes e_i \in {\cal D}$,
one has
$\displaystyle\lim_{\cal U}\lambda_i$  exists in  $\bf K$; then  setting
$\chi_{{\cal U}}(u)=
\displaystyle\lim_{\cal U}\lambda_i$, it is readily seen that
$\chi_{{\cal U}}$ is a character
of ${\cal D}$. Moreover, for all $J \in {\cal U}$ one has $\chi_{{\cal
U}}(p_J)=
\displaystyle\lim_{\cal U}1=1$, that is $J \in {\cal U}_{{\chi}_{\cal
U}}$  and
${\cal U}={\cal U}_{{\chi}_{\cal U}}$.
The theorem is proved if one considers  on
${\cal X}({\cal D})$  the weak$\mbox{ }^*$-topology 
and on $\Phi _c$ the topology
induced by the natural topology on the space of ultrafilters, which
is the weakest topology on $\Phi _c$ relative to which the mapping
$\lim : \Phi _c\to \bf K$ is continuous.
\par {\bf 2. Remark.} $(i)$. If $\bf K$ is locally compact, then  for any
bounded family
$(\lambda _i)_{i \in I} \subset \bf K$, the limit $\displaystyle \lim_{\cal
U}\lambda_i$  exists in $\bf K$. Therefore, 
$\Phi _c$ is equal to the 
entire set of all ultrafilters on $I$ and ${\cal
X}({\cal D})$ is compact,
homeomorphic to the Stone-\v Cech  compactification  $\beta (I)$ of 
the discrete topological space $I$.
\par $(ii)$. If $\bf K$ is not spherically complete and $I$ is a
small set: i.e.
the cardinal of $I$ is nonmeasurable, it is well known that the
continuous dual of
$\ell ^{\infty }(I,{\bf K})$ is equal to the space $c_0(I,{\bf K})$ of the
families converging to zero
(cf. \cite{roo} Theorem 4.21 ). Then, one can prove that ${\cal X}({\cal
D})$ is homeomorphic with $I$.
\par {\bf 3. Note.} For $\bf K$ spherically complete, not locally
compact, it is interesting to find explicit
conditions on an ultrafilter ${\cal U}$ in such a way that
$\displaystyle \lim_{\cal U}\lambda _i$
exists for any bounded family  $(\lambda _i: i \in I) \subset \bf K$. One
can try to use  Banach
limits, i.e. continuous linear forms on $\ell ^{\infty}(I,{\bf K})$ that
extend the usual continuous
linear form $lim$ defined on the subspace $c_v(I,{\bf K})$ of convergent
families. \\
Let $(J_\nu : \nu \in \Lambda )$ be a family  of subsets of $I$,
such that $J_\nu \cap J_\mu =\emptyset$ for $\nu \not= \mu$.
Putting $p_{\nu }= \displaystyle\sum_{i \in J_\nu }e'_i\otimes e_i$, one
obtains $p_{\nu }\circ p_{\mu }=\delta _{\nu, \mu }p_{\nu }$, 
for  $\nu \not= \mu $.
Hence the subalgebra with the unity ${\cal B}$
of ${\cal L}_0(E)$,  generated by $(p_{\nu }: \nu \in \Lambda )$ is equal
to ${\bf K}.id \oplus
(\displaystyle\oplus_{\nu \in \Lambda}K.p_\nu )$.
Indeed if $ u =\alpha_0 id + u_1 $ and
$ v =\beta_0 id + v_1 $  with $u_1=\displaystyle\sum_{\nu \in \Lambda}
\alpha_\nu p_\nu$ and
$v_1=\displaystyle\sum_{\nu \in \Lambda} \beta_\nu p_\nu$ (finite sums),
one has
$u \circ v=\alpha_0\beta_0 id + \alpha_0 v_1 +\beta_0 u_1 + u_1 \circ
v_1= \alpha _0\beta _0 id + \displaystyle\sum_{\nu \in \Lambda}
(\alpha _0\beta _{\nu }+ \alpha _{\nu }\beta _0 + 
\alpha _{\nu }\beta _{\nu })p_{\nu }\in {\cal B}$.                            
On the other hand, since  $u=\alpha _0 id +\displaystyle\sum_{\nu
\in \Lambda}p_{\nu }$  with
$\Gamma=\{\nu : \quad \nu \in \Lambda ; \quad \alpha_\nu \not= 0 \}$ finite
and $I=(\displaystyle\bigcup_{\nu \in \Gamma}J_\nu) \bigcup
(\displaystyle\bigcap_{\nu \in \Gamma}{J_\nu}^c)$ (a partition),  one
has $u=\alpha_0\displaystyle\sum_{i \in I}e'_i\otimes e_i +$ 
$\displaystyle\sum_{\nu \in \Gamma}\alpha_\nu \displaystyle\sum_{i \in
J_\nu}e'_i\otimes e_i $
$=\alpha_0 \displaystyle\sum_{i \in \cap _{\nu \in \Gamma}{J_\nu}^c}
e'_i\otimes e_i +$
$\displaystyle\sum_{\nu \in \Gamma}
\displaystyle\sum_{i \in J_\nu}(\alpha_0+\alpha_\nu)e'_i\otimes e_i$.
Hence $\Vert u \Vert=
\max(\vert \alpha_0 \vert, \displaystyle\max_{\nu} \vert \alpha_0
+\alpha_\nu \vert)$.
\par {\bf 4. Lemma.} {\it  
Let $ u =\alpha_0 id  + \displaystyle\sum_{\nu \in \Lambda}
\alpha_\nu p_\nu \in {\cal B}$, and
$\Lambda_0 = \Lambda \cup \{0 \}$.
Then $ \Vert u \Vert=\displaystyle\max_{\nu \in \Lambda_0}\vert
\alpha_\nu \vert$.
In other words $\{ id \}\cup \{ p_{\nu }: \nu \in \Lambda \}$ is an
orthonormal family in ${\cal L}_0(E)$.}
\par {\bf Proof.} Since $\Vert u \Vert=
\max(\vert \alpha_0 \vert, \displaystyle\max_{\nu \in \Lambda} \vert
\alpha_0 +\alpha_\nu \vert)$, and
$\displaystyle\max_{\nu} \vert \alpha _0 +\alpha _{\nu }\vert \leq
\max(\vert \alpha _0\vert,
\displaystyle\max_{\nu \in \Lambda}\vert\alpha _{\nu }\vert)$. One has
$\Vert u \Vert \leq
\displaystyle\max_{\nu \in \Lambda _0} \vert \alpha _{\nu }\vert$ .
Moreover,  $\vert \alpha _0 \vert \leq \Vert u \Vert$. Hence  for $
\nu \in \Lambda$, one
has $\vert \alpha _{\nu }\vert=\vert \alpha _{\nu }+
\alpha _0 -\alpha _0 \vert
\leq \max(\vert \alpha _{\nu }+\alpha _0\vert ,\vert \alpha _0 
\vert) \leq \Vert u \Vert$. It follows that
$\displaystyle\max_{\nu \in \Lambda _0}\vert \alpha _{\nu }\vert \leq
\Vert u \Vert $, and  Lemma 4 is proved.
\par {\bf 5. Lemma.} {\it Assume that 
$(e_i: i \in I)$ is an orthonormal basis of $E$
or $E$ is an ultrametric
Hilbert space. Then any $u \in {\cal B}$ is self-adjoint, i.e. $u^*=u$,
and $\Vert u^2 \Vert= \Vert u \Vert^2 $. }
\par {\bf Proof.} That any element of ${\cal B}$   is self-adjoint is
easy to verify.
Let $ u =\alpha _0 id  + \displaystyle\sum_{\nu \in \Lambda}
\alpha _{\nu }p_{\nu }\in {\cal B}$,  one has
$u^2={\alpha _0}^2 id  + \displaystyle\sum_{\nu \in \Lambda}
(2\alpha _0 \alpha _{\nu }+{\alpha _{\nu}}^2) p_{\nu }\in {\cal B}$.
Hence $\Vert u^2\Vert=
\max (\vert \alpha_0 \vert^2, \displaystyle\max_{\nu \in \Lambda}\vert
{\alpha_0}^2+ 2\alpha_0 \alpha_\nu +{\alpha_\nu}^2 \vert)=
(max (\vert \alpha_0 \vert , \displaystyle\max_{\nu \in
\Lambda}\vert \alpha_0+\alpha_\nu \vert))^2 =\Vert u \Vert^2$.
\par {\bf 6. Note.} In fact, Lemma 5 is true for $ u \in {\cal D}$.
Let $E$ be a free  Banach space with orthogonal basis
$(e_i: i \in I)$. Fix
$\pi \in K$ such that $0 < \vert \pi \vert < 1$. There exists  for
any $i \in I$ an integer
$n_i \in {\bf Z}$ such that $\vert \pi \vert ^{n_i+1} < \Vert e_i
\Vert \leq
\vert \pi \vert^{n_i}$. For $x=\displaystyle \sum_{i\in I}x_ie_i$, one
has $\displaystyle \lim_{i \in I}x_i\pi ^{n_i}=0$. Hence
one defines on $E$ a norm by setting $\Vert x
\Vert_{\pi }=\displaystyle \sup_{i \in I}
\vert x_i \vert \vert \pi \vert ^{n_i}$; this norm is equivalent to
$\Vert \quad  \Vert$  with
$\vert \pi \vert \Vert x \Vert _{\pi }\leq \Vert x \Vert \leq \Vert x
\Vert _{\pi }$.
Furthermore, setting for $x=\displaystyle \sum_{i \in I}x_ie_i$
and $y=\displaystyle \sum_{i \in I}y_ie_i \in E,
f_{\pi }(x,y)=\displaystyle \sum_{i \in
I}\pi ^{2n_i}x_iy_i$, one has a continuous, non degenerated,
bilinear form on  $E$ such that  $\vert f_{\pi }(x,y) \vert
\leq \Vert x \Vert _{\pi }\Vert y \Vert _{\pi }\leq \vert \pi 
\vert ^{-2} \Vert x \Vert \Vert y \Vert $.
Therefore, one obtains on $E$, a structure of ultrametric
Hilbert space $E_\pi =(E,
\Vert \quad \Vert_\pi ,f_{\pi })$. Since the norms $\Vert \quad \Vert $
and $\Vert \quad \Vert_{\pi }$
are equivalent, one  has ${\cal L}(E)={\cal L}(E_{\pi })$ and  ${\cal
L}_0(E)={\cal L}_0(E_{\pi })$.
The norms on ${\cal L}(E)$ induced by $\Vert \quad \Vert$ and  $\Vert
\quad \Vert_{\pi }$ are
equivalent with $\vert \pi \vert \Vert u \Vert_\pi \leq \Vert u \Vert
\leq \vert \pi \vert^{-1}
\Vert u \Vert _{\pi }$.
As in \S 3.4, one defines the adjoint $u^*$ of $u \in {\cal
L}(E)$ with respect to
$f_\pi$. One obtains the results stated in Theorem 3.5, that is : $u$
admit an  adjoint with respect
to $f_{\pi }$ if and only if $u \in {\cal L}_0(E)$. Furthermore
if  $u=\displaystyle \sum_{i,j}\alpha _{ij}
{e'}_j\otimes e_i \in {\cal L}_0(E)$, one has $u^*=\displaystyle \sum_{i,j}
\pi ^{n_j-n_i}\alpha _{ji}
{e'}_j\otimes e_i$, and $u$ is self-adjoint, i.e. $u^*=u$ if and only if
$\pi ^{n_i}\alpha _{ij} = \pi ^{n_j}\alpha _{ji}$, for all $i , j \in I$. 
\par {\bf 7. Note.} Let  $\pi '$ be another element of $\bf K$ such
that $0 < \vert \pi ' \vert
< 1$; let also $(m_i: i \in I) \subset {\bf Z}$ be defined by $\vert \pi'
\vert ^{m_i +1 } <
\Vert e_i \Vert \leq \vert \pi' \vert ^{m_i}$. Then the adjoint
$u^{\dagger }=
\displaystyle \sum_{i,j}{\pi '}^{m_j-m_i}\alpha _{ji}{e'}_j\otimes e_i$ 
of $u$ with respect to $f_{\pi '}$
coincides  with $u^*$ if and only if  $\pi ^{n_j-n_i}\alpha _{ji}=
{\pi '}^{m_j-m_i}\alpha _{ji}$, for each $i$ and $j \in I$. 
If this is true for all $u \in {\cal L}_0(E)$,
one has $\pi ^{n_j-n_i}={\pi '}^{m_j-m_i}$, for $i, j \in I$. Hence,
${log \vert \pi \vert \over log \vert {\pi '} \vert } = {m_j-m_i \over
n_j-n_i} = {m \over n } > 0$
and the sets  $(m_j-m_i)_{i \not= j}$ and  $(n_j-n_i)_{i \not= j}$ must
be finite.
\par If $J$ is a subset of $I$, the projector $p_J=\displaystyle\sum_{i
\in J}e'_i\otimes e_i$
is self-adjoint with respect to any bilinear symmetric form
$f_\pi$ and $\Vert p_J \Vert = 1 =\Vert p_J \Vert_\pi$.
\par {\bf 8. Lemma.} {\it Let $E$ 
be a free  Banach space with orthogonal basis
$(e_i: i \in I)$. Defining adjoint
of a continious operator with respect to $f_\pi$, one has that any $u
\in {\cal B}$ (respectively
${\cal D}$) is  self-adjoint and $\Vert u^2 \Vert =\Vert u \Vert^2$. }
\par {\bf Proof.} It is the same as in Lemma 5.
Since for any $u=\alpha _0id +\displaystyle \sum_{\nu \in
\Lambda }\alpha _{\nu }p_{\nu }
\in {\cal B}$ one has $\Vert u \Vert=
\max_{\nu \in \Lambda _0}\vert \alpha _{\nu }\vert $, 
i.e. $\{ id,p_{\nu }: \nu \in \Lambda \}$ is an
orthonormal family in ${\cal L}_0(E)$, one sees that the closure
${\cal A}=\overline{\cal B}$
of ${\cal B}$ is the subspace of ${\cal L}_0(E)$  of all elements $u$
which can be written in the unique form of summable families
$u=\alpha _0id +\displaystyle\sum_{\nu \in \Lambda}\alpha _{\nu }p_{\nu }$
with $\alpha _0,\alpha _{\nu }\in \bf K$ and
$\displaystyle\lim_{\nu }\alpha  _{\nu }=0$.
It is readily seen that ${\cal A}$ is a closed  unitary
subalgebra of ${\cal L}_0(E)$, contained in ${\cal D}$,
such that any  element  $u$  of  ${\cal A}$ is self-adjoint. Moreover
for the pointwise convergence, one has 
$u=\alpha _0\displaystyle\sum_{i \in
\cap {J_{\nu }}^c}e'_i\otimes e_i+
\displaystyle\sum_{\nu \in \Lambda}\displaystyle\sum_{i \in J_\nu }
(\alpha _0+\alpha _{\nu })e'_i\otimes e_i$. Hence, if $\bigcap_{\nu \in
\Lambda}{J_ {\nu }}^c =\emptyset $, one
has $u= \displaystyle\sum_{\nu \in \Lambda}\displaystyle\sum_{i \in
J_\nu}(\alpha_0+\alpha_\nu)
e'_i\otimes e_i$ and $id = \displaystyle\sum_{\nu \in \Lambda}p_{\nu }$.
\par {\bf 9. Example.} If $\Lambda =I$ and $J_i= \{i\}$ for each
$i \in I$, one has ${\cal A}=
\{ \alpha _0id +\displaystyle\sum_{\i \in I}\alpha _i e'_i\otimes e_i :
\quad \alpha _i \in {\bf K}, 
\quad \displaystyle\lim_{i \in I}\alpha _i =0 \}$. As an element of
${\cal D}$ any $u \in {\cal A}$ is
in the form $u= \displaystyle\sum_{i \in I}a_ie'_i\otimes e_i$ with
$\displaystyle\lim_{i \in I}a_i = \alpha _0$ exists in $\bf K$.
\par {\bf 10. Proposition.} {\it $(i)$. Any element $u$ of the  
Banach algebra ${\cal A}$ with the unit element is
self-adjoint with respect to
any bilinear symmetric form  $f_{\pi }$ and  $\Vert u^2 \Vert= \Vert u
\Vert ^2$.
\par $(ii)$. The spectrum ${\cal X}({\cal A})=Hom_{alg}({\cal A},
{\bf K})$ of ${\cal A}$, equipped with
the weak*-topology, is homeomorphic to the Alexandroff
compactification of the  discrete space $\Lambda$. }
\par {\bf Proof.} The first part is an easy consequence of Lemma 8.
Let $\chi \in {\cal X}({\cal A})$, then $\chi $ is a continuous
linear form with norm
$\Vert \chi \Vert = 1$. Furthermore, one has $\chi (id)=1$ and
$\chi (p_{\nu }p_{\mu })=
\chi (p_{\nu }) \chi (p_{\mu })=\delta_{\nu ,\mu}
\chi (p_{\nu })$, for $\nu , \mu \in \Lambda$. It follows
that for any $\nu \in \Lambda$, one has $\chi (p_{\nu })=1$ or
$\chi (p_{\nu })=0$. Hence : $(a)$ there
exists $\nu \in \Lambda$ such that $\chi (p_{\nu }) =1$ and 
$\chi (p_{\mu })=0$ for $\mu \not= \nu$,  or
\quad $(b)$ \quad $\chi(p_\nu)=0$ for all $\nu \in \Lambda $.
In the case $(a)$ one puts $\chi =\chi _{\nu }$ and in the case $(b),$
$\chi =\chi _0$. One verifies that for $u= \alpha _0id
+ \displaystyle \sum_{\nu \in \Lambda}\alpha _{\nu }p_{\nu }
\in {\cal A}$, one has $\chi _0(u)=\alpha  _0$
and $\chi _{\nu }(u)=\alpha _0 +\alpha _{\nu }$, 
$\nu \in \Lambda $. It  follows
that ${\cal X}({\cal A})=
\{ \chi _0  ,\quad  \chi _{\nu }: \quad  \nu \in  \Lambda \}$ 
and  ${\cal X}({\cal A})$ is in a
bijective correspondence with the set $\Lambda _0=\Lambda \cup \{0\}$.
Let $W(\chi ;\quad \varepsilon ,\quad u_1,......,u_n) =
\{ \eta: \quad \eta  \in {\cal X}({\cal D});
\quad \vert \chi (u_j) - \eta (u_j) \vert < \varepsilon,  
\quad u_j \in {\cal A},\quad 1 \leq j \leq n \} $
be a fundamental neighborhood  of $\chi \in  {\cal X}({\cal A})$ for the
weak*-topology. Since for
$u_j=\alpha _{0j}id+\displaystyle \sum_{\mu \in \Lambda}\alpha _{\mu
j}p_\mu \in {\cal A}$, one has
$\displaystyle \lim_{\mu \in \Lambda }\alpha _{\mu j} = 0$,  there exists
a finite subset
$\Gamma _{\varepsilon }$  of
$\Lambda $,  such that for any $\mu \not\in  \Gamma _{\varepsilon }$, one
has $\vert \alpha _{\mu j} \vert < \varepsilon $ for each $1 \leq j \leq n$.
If $\chi =\chi _{\nu }, \nu \in \Lambda $,  one has for $1 \leq  j
\leq n, \mu \in \Lambda ,
\chi _{\nu }(u_j) - \chi _{\mu }(u_j) =\alpha _{\nu j} - \alpha _{\mu j}$.
Choosing $(u_j: 1 \leq  j \leq n)$ such that $\varepsilon _{\nu }=
\displaystyle \min_{1 \leq j \leq n}\vert \alpha _{\nu j} \vert > 0$,
there exists $\Gamma _{\nu }\subset \Lambda , \Gamma _{\nu }$
finite  such that $\vert \alpha _{\mu j} \vert < \varepsilon _{\nu }$ 
for $1 \leq j \leq n$ and for all $\mu \not\in \Gamma _{\nu }$. 
Hence $\vert \alpha _{\mu j} \vert < \vert
\alpha _{\nu j} \vert $
and $\vert \alpha _{\nu j}-\alpha _{\mu j} \vert = \vert \alpha _{\nu j}
\vert \geq \varepsilon _{\nu }$,
for $1 \leq j \leq n$ and $\mu \not\in \Gamma _{\nu }$. Therefore, if
$\varepsilon < \varepsilon _{\nu }$, one
has  $W(\chi _{\nu };\quad \varepsilon , \quad u_1,......,u_n) = \{
\chi _{\nu } \}$, that is
$\{ \chi _{\nu } \}$ is open in ${\cal X}({\cal A})$. Hence  
$\{ \chi _{\nu }: \nu \in \Lambda \}$ 
is a discrete subset of ${\cal X}({\cal A})$.
On the other hand, if $\chi =\chi _0$, one has $\chi _0 (u_j) -
\chi _\mu (u_j) =
-\alpha _{\mu j}$. Hence for $\varepsilon > 0$, there exists a finite
subset $\Gamma _{\varepsilon }$ of
$\Lambda $ such that for $\mu \not\in \Gamma _{\varepsilon }$, one has
$ \vert \chi _0 (u_j) - \chi  _{\mu }(u_j) \vert =
\vert \alpha _{\mu j} \vert < \varepsilon $ for each
$1 \leq j \leq n $. In other words, $W(\chi _0 ; \quad \varepsilon ,
\quad u_1,......,u_n) =
\{ \chi _{\mu }: \mu \not\in  \Gamma _{\varepsilon } \}$. Furthermore
$\chi _0 =
\displaystyle\lim_{\mu \in \Lambda}\chi _{\mu }$ in  ${\cal X}({\cal A})$
for the weak*-topology.
It follows that ${\cal X}({\cal A})$ is weak*-compact.
Consider on $\Lambda _0=\Lambda \cup \{  0 \}$ the topology such
that $\Lambda $ is a discrete
subset of $\Lambda _0$ and the neighborhoods of $0$ are $W_{\Gamma }(0) =
\Lambda _0 \setminus \Gamma $, 
where $\Gamma \subset \Lambda $ is finite. It
becomes  clear that
$\Lambda _0$ is homeomorphic to the Alexandroff compactification of the
discrete space  $\Lambda $.
Identifying ${\cal X}({\cal A})$ with $\Lambda _0$ one concludes the
proof of the proposition.
\par {\bf 11.} Let 
${\cal C}({\cal X}({\cal A}),{\bf K})$ be the $\bf K$-Banach algebra
of the continuous functions
$f$ on the compact space ${\cal X}({\cal A})$ with values in $\bf K$. It is
readily seen that
$f \in {\cal C}({\cal X}({\cal A}),{\bf K})$ is defined by the family
$(f(\chi _{\nu }): \nu \in
\Lambda _0) \subset \bf K$ such that $\displaystyle\lim_{\nu \in
\Lambda}f(\chi _{\nu }) = f(\chi _0)$.
Hence ${\cal C}({\cal X}({\cal A}),{\bf K})$ is isometrically isomorphic
to the algebra
$c_v(\Lambda _0,{\bf K})= \{ a: \quad a=(a_{\nu }: \nu \in \Lambda _0) 
\subset {\bf K }; \quad
\displaystyle\lim_{\nu \in \Lambda } a_\nu =a_0 \}$ : on  $c_v(\Lambda _0,
{\bf K})$, one considers the
usual multiplication defined pointwise and the norm 
$\Vert (a_{\nu }: \nu \in \Lambda _0) \Vert =
\displaystyle\sup_{\nu \in \Lambda _0}\vert a_{\nu }\vert $
\par {\bf 11.1. Corollary.} {\it  The Banach algebra 
${\cal A}$ with the unit element is isometrically isomorphic to the algebra
$c_v(\Lambda _0,{\bf K})= \{ a: \quad a=
(a_\nu)_{\nu \in \Lambda_0} \subset {\bf K}; 
\quad \displaystyle\lim_{\nu \in \Lambda } a_{\nu }=a_0 \}$   }.
\par {\bf Proof.} Let ${\cal G}: 
{\cal A}  \longrightarrow   {\cal C}({\cal X}({\cal
A}),{\bf K})$  be the Gelfand
transform:  ${\cal G}(u)(\chi ) = \chi (u)$. As usual, ${\cal G}$ is
continuous. Since for $u = \alpha _0id
+ \displaystyle\sum_{\nu \in \Lambda }\alpha _{\nu }p_{\nu }\in {\cal A}$,
one has $\chi _0 (u)=\alpha _0$
and $\chi _{\nu }(u) = \alpha _0 + \alpha _{\nu }, \nu \in \Lambda $, one
obtains $\Vert u \Vert =
\max (\vert \chi _0 (u) \vert , \displaystyle\sup_{\nu \in \Lambda}\vert
\chi _{\nu }(u) \vert )=
\displaystyle\sup_{\chi \in {\cal X}({\cal A})} \vert \chi (u) \vert $.
Hence, $\Vert {\cal G}(u) \Vert = \Vert u \Vert$.
Furthermore, ${\cal G}(id)(u) =1$, i.e. ${\cal G}(id)=f_0$ the
constant function equal to $1$.
On the other hand, for $\nu \in \Lambda  $, ${\cal G}(p_{\nu })(\chi )=1$
if $\chi =\chi _{\nu }$ and $0$ otherwise. 
Hence, setting for $\nu \in \Lambda $,
$f_\nu : {\cal X}({\cal A}) \longrightarrow \bf K$ such that
$f_{\nu }(\chi _{\mu })=\delta _{\nu,\mu}, \mu
\in \Lambda $, one has ${\cal G}(p_{\nu })=f_{\nu }$.
Let $u=\alpha _0id + \displaystyle\sum_{\nu \in \Lambda}\alpha_{\nu }
p_{\nu }\in {\cal A}$, one has
${\cal G}(u)=\alpha _0f_0+\displaystyle\sum_{\nu \in \Lambda}\alpha _{\nu }
f_{\nu }$. Since any $f \in {\cal C}({\cal X}({\cal A}),{\bf K})$ 
can be written in the unique convergent sum
$f= f(\chi _0)f_0+\displaystyle\sum_{\nu \in
\Lambda }(f(\chi _{\nu })-f(\chi _0))f_{\nu }$ with
$\displaystyle\lim_{\nu \in \Lambda}(f(\chi _{\nu })-f(\chi _0)) =0$, one
has $f={\cal G}(u)$ with
$u=f(\chi _0)id + \displaystyle\sum_{\nu \in
\Lambda }(f(\chi _{\nu })-f(\chi _0))
p_{\nu }$. Hence, ${\cal G}$ is surjective. 
Together with $\Vert {\cal G}(u) \Vert = \Vert
u \Vert$, the Corollary is proved.
\section{Spectral integration.}
\par {\bf 1.} Suppose that $X$ and $Y$ are Banach spaces 
over a topologically complete non-Archimedean field $\bf K$ 
with a non-trivial valuation and ${\cal L}(X,Y)$ denotes the Banach space
of bounded linear operators $E: X\to Y$ supplied with the operator norm:
$\| E \| :=\sup_{0\ne x\in X} \| Ex\| _Y/\| x\| _X$. For $X=Y$
we denote ${\cal L}(X,Y)$ simply by ${\cal L}(X)$.
Let $X$ and $Y$ be isomorphic with the Banach spaces
$c_0(\alpha ,{\bf K})$ and $c_0(\beta ,{\bf K})$ 
and let they be supplied with 
the standard orthonormal bases $\{ e_j: j\in \alpha \} $
in $X$ and $\{ q_j: j\in \beta \} $ in $Y$ respectively, 
where 
$c_0(\alpha ,{\bf K}):=\{ x=(x^j: j\in \alpha)|$ $x^j\in {\bf K}, $
$\mbox{ such that for each }\epsilon >0$ $\mbox{ a set }
\{ j: |x^j|>\epsilon \} \mbox{ is finite } \} $ with a norm
$\| x\|:=\sup_j|x^j|_{\bf K}$, $\alpha $ and $\beta $ are ordinals
(it is convenient due to Kuratowski-Zorn Lemma).
Then each operator $E\in {\cal L}(X,Y)$
has its matrix realisation $E_{j,k}:=q_k^*Ee_j$, 
which may be infinite, where
$q_k^*\in Y^*$ is a continuous
$\bf K$-linear functional $q_k^*: Y\to \bf K$ corresponding to $q_k$
under the natural embedding $Y\hookrightarrow Y^*$
associated with the chosen basis, $Y^*$ is a topologically
conjugated or dual space of $\bf K$-linear functionals on $Y$,
$q_k^*(q_l)=\delta ^k_l$.
Therefore, to each $E\in {\cal L}(X,Y)$ there corresponds an adjoint operator
$E^*\in {\cal L}(Y^*,X^*)$. 
By a transposed operator $E^t$ we mean a restriction $E^*|_Y$,
where $Y$ is embedded into $Y^*$
such that $E^t_{j,k}=E_{k,j}$ for each $j\in \alpha $
and $k\in \beta $. 
\par This means, that if $X=Y$ and $E^t=E$, then $E$ is called a symmetric operator.
For $X=Y=c_0(\alpha ,{\bf K})$ there is an inclusion
$E^*(Y)\subset X^*$. Since $X^*=l^{\infty }(\alpha ,{\bf K})$,
then $\| x \|_X=\| x\|_{X^*}$
for each $x\in X$. 
Since $\| E\| =\sup_{j,k} |E_{j,k}|$, then
$\| E\| =\| E^*\|$ and $\| E\| =\| E^t\| $.
If $A$, $E\in {\cal L}(X)$ and $E=A^t$, then
$A$ and $E$  belong to the closed subalgebra
${\cal L}_0(X)$ (see \S 2).
\par {\bf 1.1.} Now let $A$ be an abstract Banach algebra over a field 
$\bf K$, which is complete relative to a norm $\| * \| _{\bf K}$
in it. We say that $A$ is with an operation of transposition $a\mapsto a^t$
for each $a\in A$ if the following conditions $(\alpha - \delta )$ are
satisfied:
\par $(\alpha )$ $(a+b)^t=a^t+b^t;$
\par $(\beta )$ $(\lambda a)^t=\lambda a^t;$
\par $(\gamma )$ $(ab)^t=b^ta^t;$
\par $(\delta )$ $(a^t)^t=a^{tt}=a$
for each $a, b\in A$ and each $\lambda \in \bf K$.
\par Let $A$ be an algebra over $\bf K$, 
which satisfies the following conditions $(i-iii):$
\par $(i)$ $A$ is a Banach algebra
\par $(ii)$ with
the operation of transposition $a\mapsto a^t$,
\par $(iii)$ $\| a^ta \|=\| a\| ^2$ for each $a\in A$
while evaluation of norms. 
\par Then such algebra is called an $E$-algebra.
\par Without condition $(iii)$ it is called the $T$-algebra. 
If instead of $(iii)$ it is satisfied the following condition:
\par $(iv)$ $\| a^ta \| =\| a^2 \| $ then $A$ is called a $S$-algebra.
\par For each $E$-algebra we have $\| a \| ^2=\| aa^t\| 
\le \| a \| \| a^t \| $, hence $\| a \| \le \| a^t \| $ and also
$\| a^t \| \le \| (a^t)^t \| = \| a \| $, consequently, 
$\| a \| = \| a^t \| $.
\par {\bf 1.2.} Evidently, ${\cal L}_0(X)$ is the $T$-algebra.
Each $C$-algebra is at the same time the $E$-algebra (see also
\S 5.6), since for each singleton $x\in X$ a closed subalgebra
$C(\{ x \} , {\bf K} )$ is isomorphic with $\bf K$ and the 
restriction of transposition on $C(\{ x \} , {\bf K} )$ 
gives $f^t(x)=f(x)$ for each $f\in C_{\infty }(X,{\bf K})$.
\par Let $A_1$,..., $A_n$ be linear operators. Then 
the equation $[\sum_{j=1}^n\lambda _jA_j, \sum_{k=1}^n\mu _kA_k ]
=\sum_{j<k}(\lambda _j\mu _k -\lambda _k \mu _j)(A_jA_k-A_kA_j)=0$
for each $\lambda _j$ and $\mu _k \in \bf K$ is equivalent to
$[A_j,A_k]=0$ for each $j<k$.
In view of \S \S IV.6,7, VII.7, VIII.2 \cite{gant} for each $n\in \bf N$
there are pairwise commuting matrices of the size $m\times m$
such that for sufficiently large $m>n$ they in addition can be found 
non-diagonal (non-reducible to diagonal form by
transformations, $U_jA_jU^{-1}_j$, where $U_j$ are invertible matrices),
since in view of Theorem VIII.7.2 \cite{gant} a number of linearly 
independent matrices, which commute with the given matrix 
$A$ is defined by the following formula: $N=n_1+3n_2+...+(2t-1)n_t$, 
where $n_1$, $n_2$,...,$n_t$ are degrees of non-constant invariant
polynomials $i_1(\lambda )$, $i_2(\lambda )$,...,$i_t(\lambda )$
and $n=n_1+n_2+...+n_t$ is a size of the square $n\times n$ matrix $A$.
This can be done by suitable choices
of Jordan forms of matrices over $\bf K$. From this it follows, that
in ${\cal L}_0(c_0(\omega _0,{\bf K}))$ for each $n\in \bf N$
there always exist $n$ pairwise commuting operators such that
they are not reducible to the diagonal form by adjoint transformations
$U_jA_jU^{-1}_j$. This produces examples of $T$-algebras.
When in the (finite case ) 
Jordan form $| \lambda _{j,k} | >1 $ and $\| I - U_j \| <1$
for each $j=1,...,n$, where $\lambda _{j,k}$ are diagonal elements 
of the Jordan normal forms of $A_j$, then each $A_j$ satisfy condition
$(iii)$ together with $A_j^t$. We take the case
$U_j=U_1=:U$ for each $j$.
Let in addition $A_j$ be pairwise commutative and
have block forms $(v)$ $A_j=diag(A_{j,1},...,A_{j,n})$,
$A_{j,l}$ are $m_l\times m_l$ square matrices and $A_{j,l}=0$ for
each $l\ne j$ and $A_{j,j}\ne 0$. Consider their transposed matrices also, 
then the linear span of all their products $B_{j_1}...B_{j_s}$ 
produces commutative $E$-algebra, which is generated by non-diagonalizable
over $\bf C_p$ matrices, where $B_j$ is equal either to $A_j$ or to
$A_j^t$, since $|\lambda _{j,k} |^a>|b\lambda _{j,k}|^c$ for each $a>c>0$
and $b\in \bf Z$. As it follows from \S VIII.7.2 \cite{gant}
these $A_j$ can be chosen such that $sp_{\bf K}\{ I, A_j, A_j^2,...,A_j^m \} $
does not contain $A_l$ for each $l\ne j$.
\par The construction of $S$-algebras can be done analogously
and more lightly, since Condition $(iii)$ is replaced by
Condition $(iv)$.
\par Let $A$ be an $E$-algebra with a $\bf K$-linear 
isometry $Y: A\to A$ such that $Y(ab)=Y(b)Y(a)$ and
$Y^t(a)=Y(a^t)$ for each $a, b\in A$. Then $Y(A)$ is the $E$-algebra.
\par There are general constructions of Banach algebras also.
In particular we can take a free Banach algebra $A$ 
generated by a set $J$. This means, that $A$ is a completion of
$sp_{\bf K} \{ a_1...a_n:$ $a_1,..,a_n\in J,$ ${n\in \bf N} \} $
with the definite order of letters $a_1,...,a_n$ in each word
$w=a_1a_2...a_n$, when neighbouring elements $a_j$ 
and $a_{j+1}$ are distinct in $J$. There exists a norm on 
$sp_{\bf K} \{ a_1...a_n:$ $a_1,...,a_n\in J,$ ${n\in \bf N} \} $
such that $\| ab\| \le \| a\|$ $\| b \| $ for each $a, b\in A$.
For example, $\| w\| =1$ for each word $w=a_1...a_n$ 
with $a_1,...,a_n\in J$,
$\| c_1w_1+...+c_mw_m \| =\max_{1\le j\le m} |c_j|_{\bf K}$
for different words $w_1,..., w_m$ with $c_j\in \bf K$
for each $j=1,...,m$.
Then for $Y: A\to A$ preserving a closed ideal $V$ we can consider 
the quotient mapping $\bar Y: \bar A\to \bar A$, where
$\bar A=A/V$ and $Y$ on $A$ is defined by $Y$ on $J$ due to
the continuous extension.
\par Another example is the following. For a subset $J$ of symmetric
(that is, $a^t=a$ for each $a\in J\subset {\cal L}_0(X)$) 
pairwise commuting elements 
(that is, $ab=ba$ for each $a$ and
$b\in J\subset {\cal L}_0(X)$) let $A:=cl\mbox{ }
(sp_{\bf K}\{ \prod_{i=1}^ma_i^{n_i}: 0\le n_i\in {\bf Z}, m\in {\bf N},
a_i\in J \} )$,
where $a^0:=I$ is the unit operator on $X$.
Then such $A$ is the $T$-algebra. Since $\| a^t a\|=\| a^2 \| $
for each $a\in A$, then it is the $S$-algebra.
\par Another example of an $E$-algebra is the algebra of diagonal
operators in ${\cal L}_0(X)$. Then each $E$-algebra is certainly
a $S$-algebra and each $S$-algebra is a $T$-algebra
(see also Lemmas 4.5, 4.8 and Proposition 4.10.(i)).
Above were constructed more interesting examples of $E$-algebras and
$S$-algebras. In general diagonal form of an algebra is unnecessary for
the spectral theory. Moreover, there are well-known theorems,
when Lie algebras (in particular of finite square $m\times m$ 
matrices over $\bf C_p$)
can be reduced simultaneously to the upper triangular form
by one transformation $UA_jU^{-1}$ (see K. Iwasawa Theorem 
4.7.2 \cite{gogro}).
There are also cases, when it may be over $\bf K$.
Using limits such cases can be spread on subalgebras of
${\cal L}_0(X)$.
\par It will be shown below that for the spectral integration
it is sufficient to consider $C$-algebras.
\par {\bf 2.} Let $A$ and $B$ be two $E$-algebras over the same field
$\bf K$,
an algebraic homomorphism $T: A\to B$ is called a $t$-representation of
$A$ in $B$, if  $T_{a^t}=(T_a)^t$  for each $a\in A$.
The reducing ideal $\sf \Upsilon $ of $A$ is defined as the intersection of the kernels
of all $t$-representations of $A$. $\sf \Upsilon $ is also called the $t$-radical.
If ${\sf \Upsilon }=0$, then $A$ is called reduced (or  $t$-simple).
\par Let $\| a\|_t:=\sup_{T\in \Psi } \| Ta\| $ for a reduced algebra $A$, 
where $\Psi :=\Psi _A$ denotes  the family of all $t$-representations of $A$.
Since $A$ is reduced, then $\| a\|_t\ne 0$ for each $A\ni a\ne 0$.
Such $\| *\| _t$ is called a $E$-norm of $A$.
\par The $E$-algebra obtained by completing $A/{\sf \Upsilon }$ by its
$E$-norm is called the $E$-completion of $A$ and is denoted by
$A_t$.
Denote by $\pi : A\to A/{\sf \Upsilon }$ the natural $t$-homomorphism
of $A$ into $A_t$ such that $\pi (a)=a+\sf \Upsilon $ for each $a\in A$.
Then the map $T\mapsto T'=T\circ \pi $  is a bijective correspondence
between the set of all $t$-representations  $T$ of $A_t$ and the set
of all $t$-representations $T'$ of $A$.
This correspondence preserves closed stable subspaces, non-degeneracy,
bounded intertwining operators, isometric equivalence and Banach direct sums.
\par {\bf 3.} Let $A$ be a commutative Banach $T$-algebra and $A^+$ denotes
the Gelfand space of $A$, that is, $A^+=Sp(A)$, where $Sp(A)$ was defined in
Ch. 6 \cite{roo}: it is the set of all nonzero algebra homomorphisms
$\phi : A\to \bf K$ topologized as the subset of ${\bf K}^A$.
Every $x\in A$ induces a function $G_x: Sp (A)\to \bf K$ by
$G_x(\phi ):=\phi (x)$, where $\phi \in Sp(A)$, $G_x$ is called
the Gelfand transform of $x$, $G$ is called the Gelfand transfromation.
Then it is defined the spectral norm $\| x \| _{sp}:=\sup_{\phi \in Sp (A)} 
|G_x(\phi )|$ of $x\in A$. 
If $Sp (A)=\emptyset $, then $\| x \|_{sp}:=0$ for each $x\in A$.
We denote by $\hat A$ the closed subset of
$A^+$ consisting of those $\phi \in A^+$ for which $\phi (a^t)=
\phi (a)$ for each $a\in A$, $\phi \in A$ is called symmetric, if
$\phi \in \hat A$. Let $C_{\infty }(\hat A,{\bf K})$ 
be the same space as in
\cite{roo}. For a locally compact $E$ the space $C_{\infty }(E,{\bf K})$
is a subspace of the space $BUC(E,{\bf K})$
of bounded uniformly continuous functions $f: E\to \bf K$
such that for each $\epsilon >0$ there exists a compact subset
$V\subset  E$ for which $|f(x)|<\epsilon $ for each $x\in E\setminus V$.
When $E$ is not locally compact and has an embedding
into $B({\bf K},0,1)^{\gamma }$ such that $E \cup \{ x_0 \}=cl(E)$
we put $C_{\infty }(E,{\bf K}):=\{ f\in C(E,{\bf K}): 
\lim_{x\to x_0}f(x)=0 \}$,
where $B(X,x,r):=\{ y\in X: d(x,y)\le r \} $ is a ball in the metric space
$(X,d)$, $cl(E)$ is taken in $B({\bf K},0,1)^{\gamma }$, 
$\gamma $ is an ordinal, $x_0\in B({\bf K},0,1)^{\gamma }$.
\par {\bf 3.1. Definition} (see also Ch. 6 in \cite{roo}). A commutative 
Banach algebra $A$ is called
a $C$-algebra if it is isomorphic with $C_{\infty }(X,{\bf K})$
for a locally compact Hausdorff hereditarily
disconnected space $X$, where $f+g$ and $fg$
are defined pointwise for each $f, g\in C_{\infty }(X,{\bf K})$.
\par {\bf 3.2. Proposition.}  {\it  The reducing ideal $\sf \Upsilon $ of $A$ 
consists of those $a\in A$ such that $\hat a(\phi )=0$ for each
$\phi \in \hat A$. The equation 
\par $(i)$ $F(\pi (a))=\hat a|_{\hat A}$ \\
determines an isometric $t$-isomorphism $F$ of $A_t$  onto $C_{\infty }
(\hat A,{\bf K})$, when $\bf K$ is a locally compact field.}
\par {\bf Proof.}  We have $\sup_{\phi \in \hat A} |\hat a (\phi )|
\le \| \pi (a)\| _t$ for each $a\in A$.  Then we take a $t$-representation
$T$ of $A$ and $B:=cl_{\| \mbox{ }\| }range (T)$, hence $B$ is a commutative
$E$-algebra. To finish the proof of Proposition 3 we need the following.
\par {\bf 4.} {\bf Proposition.} {\it Let $A$ be a commutative $E$-algebra 
as in \S 1: $A=cl\mbox{ }(sp_{\bf K}\{ \prod_{i=1}^ma_i^{n_i}: 0\le n_i
\in {\bf Z},
a_i\in J, a^t=a, m\in {\bf N} \} )$, 
then $\hat A=A^+$. Furthemore, the Gelfand transform
map $a\mapsto \hat a$ is an isometric isomorphism of $A$ onto
$C_{\infty }(\hat A,{\bf K})$, when $\bf K$ is a locally compact field.}
\par {\bf Proof.}  In view of Corollaries  6.13, 6.14 and 6.17
\cite{roo} it is sufficient to show that $A^+=\hat A$, since
$A$ is isomorphic with $C_{\infty }(Sp(A),{\bf K})$, where $Sp(A)=\hat A$.
If $a^t=a\in A$ and $\phi \in A^+$, then $\phi (a^t)=\phi (a)\in \bf K$.
If ${\bf 1}\notin A$, then $\phi $ extends to a $t$-homomorphism
of  the $S$-algebra $A_1$ obtained by adjoining the unit $\bf 1$ to $A$,
since it is possible to consider $X\oplus {\bf K}$
(see, for example, about adjoining of $\bf 1$ in \S VI.3.10 \cite{fell}
and Ch. 6 \cite{roo}). 
Since $\bf K$ is locally compact, $A_1$ is isomorphic with $C(\alpha Y,{\bf K})$,
where $\alpha Y=Sp(A)\cup \{ 0 \} $ is a one-point (Alexandroff) compactification
of $sp(A)$ (see Observation 6.2 in \cite{roo}). Indeed, $\| a\| =\sup_{\chi
\in A^+} |\chi (a) |$.
Let $\phi (a)=r\in \bf K$,
$b:=a+z\bf 1$, where $z\in \bf K$. From $\| \phi \| \le 1$
it follows, that $| \phi (b^tb)|_p=|(r+z)^2|_p
\le \| b^2\| $ and $\| b^2 \| =\| a^ta+za^t+za+z^2{\bf 1} \| \le 
\max (\| a^2\|, |z|_p\mbox{ }\| a\|,
|z|_p^2)$. Then
there exists $0<\epsilon <(\| a^2\| )^{1/2}$ such that for each 
$|z|_p<\epsilon $: $| \phi (b^tb)|_p\le \| a^2\| $, consequently,
$\phi $ has the continuos extension on $A_1$.
\par If $\phi \in A^+$ and $a=b+c\in A$ with $b^t=b$ and $c^t=-c\in A$, 
then $\phi (a^t)=\phi (b)-\phi (c)$.
If $\phi (a^t)=-\phi (a)$ for each $a\in A$, then $\phi (b)=0$ for each
$b=b^t\in A$. If $\phi (a^t)=\phi (a)$ for each $a\in A$, then $\phi (c)=0$
for each $c^t=-c\in A$. The operation of transposition $a\mapsto a^t$
is continuous in $A$. Let $\phi \in A^+$ and $\phi \ne 0$.
Therefore, for each $\phi (a)\ne 0$ we have
$\phi (a^t)\ne 0$, since $a^{tt}=a$. Whence  $\phi (a^t)=\lambda _{\phi }
\phi (a)$ for each $a\in A$ such that $\phi (a)\ne 0$, since
$coker_{\bf K} \phi $ is one-dimensional, where 
$0\ne \lambda _{\phi }\in \bf K$. We have $\phi ((a^ta)^n)=\lambda _{\phi }^n
\phi (a)^{2n}$. Since $\| a^ta\| =\| a^2\| $, $a^{tt}=a$ and $\phi $
is the continuous multiplicative linear functional, then $|\lambda _{\phi }|_p
=1$. On the other hand, $\lambda _{\phi }\phi (ab)=\phi (a^tb^t)=\phi (a^t)
\phi (b^t)=\lambda _{\phi }^2\phi (a)\phi (b)=\lambda _{\phi }\phi (ab)$,
hence $\lambda _{\phi }=1$, where there are $a$ and $b\in A$
such that $\phi (ab)\ne 0$. 
Therefore, $\phi ^t=\phi $, where
$\phi ^t(a):=\phi (a^t)$ for each $a\in A$.
Consequently, $\hat A=A^+$.
\par {\bf Continuation of the proof of Proposition 3.2.} 
In view of \S 2 there exists $\psi \in \hat B$ 
such that $|\psi (T_a)|=\| T_a\| $, since
$\psi ': a\mapsto \psi (T_a)\in \hat A$, hence $\| T_a\|= |\psi '(a)|\le
\sup_{\phi \in \hat A}|\hat a(\phi )|$, consequently,
$\| \pi (a) \| _t\le \sup_{\phi \in \hat A}|\hat a (\phi )|$. Therefore,
$\| \pi (a)\| _t =\sup_{\phi \in \hat A}|\hat a(\phi )|$, hence the map
$F$ defines the isometric $t$-isomorphism of $A_t$ 
into $C_{\infty }(Sp(A),{\bf K})$.
The range of $F$ is a $T$-subalgebra of $C_{\infty }(\hat A,{\bf K})$, 
which automatically separates points of $\hat A$, consequently, by the 
Kaplansky theorem $cl\mbox{ }range (F)=C_{\infty }(\hat A,{\bf K})$
(see \S A.4 \cite{sch1}).
\par {\bf 5.} Let $H=c_0(\alpha ,{\bf K})$,
where $\bf K$ is a topologically complete field. 
A strong operator topology in ${\cal L}(H,Y)$ (see \S 1)
is given by a base $V_{\epsilon ;E;x_1,...,x_n}:=\{ Z\in {\cal L}(H,Y):
\sup_{1\le j\le n} \| (E-Z)x_j \|_Y <\epsilon \}$, where 
$0<\epsilon $, $E\in {\cal L}(H,Y)$, $x_j\in H$; $j=1,...,n;$ $n\in \bf N$.
Let also $X$ be a topological
space with the small inductive dimension $ind (X)=0$.
An $H$-projection-valued measure on an algebra
$\sf L$ of subsets of $X$ 
is a function $P$ on $\sf L$ assigning to each $A\in \sf L$
a projection $P(A)$ on $H$ and satisfying the following conditions:
\par $(i)$ $P(X)={\bf 1}_H$,
\par $(ii)$ for each sequence $\{ A_n : n=1,...,k \} $
of pairwise disjoint sets in $\sf L$ there 
are projections $P(A_n)$ such that $P(A_n)P(A_l)=0$
for each $n\ne l$
and $P(\bigcup_{n=1}^kA_n)=\sum_{n=1}^kP(A_n)$,
\par $(iii)$ if ${\sf A}\subset \sf L$ is shrinking and
$\cap {\sf A}=\emptyset $, then $\lim_{A\in \sf A}P(A)=0$,
where the convergence on the right hand side is unconditional in the
strong operator topology and the sum is equal to the projection
onto the closed linear span over $\bf K$
of $\{ range(P(A_n)): n=1,...,k \} $
such that $P(\emptyset )=0$, $k\in \bf N$.
\par If $\eta \in H^*$ and $\xi \in H$, then $A\mapsto \eta (P(A)\xi )$
is a $\bf K$-valued measure on $\sf L$.
The case of a $\sigma $-algebra $\sf L$ and of
$k=\infty $ in $(ii)$ is unnecessary for the subsequent consideration
and it will not be used, but it may be considered as a particular case.
The $\sigma $-additive case leads to the restriction that
each measure $\eta (P(A)\xi )$ is atomic, when $\bf K$ is spherically 
complete (see Chapter 7 in \cite{roo}).
\par Then by the definition $P(A)\le P(B)$ if and only if $range (P(A))
\subset range P(B)$. There are many projection operators on $H$,
but for $P$ there is chosen some such fixed system.
\par A subset $A\subset X$ is called $P$-null 
if there exists $B\in \sf L$ such that $A\subset B$ and $P(B)=0$,
$A$ is called $P$-measurable if $A\bigtriangleup B$ is $P$-null, 
where $A\bigtriangleup B:=(A\setminus B )\cup (B\setminus A)$.
A function $f: X\to \bf K$ is called $P$-measurable, if
$f^{-1}(D)$ is $P$-measurable for each $D$ in a field
$Bco ({\bf K})$ of clopen subsets of $\bf K$.
It is essentially bounded, if there exists $k>0$ such that
$\{ x: |f(x)|>k \} $ is $P$-null, $\| f\| _{\infty }$
is by the definition the infimum of such $k$. Then ${\sf F}:=
sp_{\bf K} \{ Ch_B: B\in {\sf L} \} $ is called the 
space of simple functions,
where $Ch_B$ denotes the characteristic function of $B$.  The 
completion of $\sf F$ relative to $\| *\| _{\infty }$ is the Banach algebra
$L_{\infty }(P)$ under the pointwise multiplication.
\par There exists a unique linear mapping
${\sf I}: {\sf F}\to {\cal L}(H)$ by the following formula:
\par $(iv)$ ${\sf I}(\sum_{i=1}^n\lambda _i Ch_{B_i})=\sum_{i=1}^n
\lambda _iP(B_i)$, where $n\in \bf N$, 
$B_i\in \sf L$, $\lambda _i\in \bf K$.
Since
\par $(v)$ $\| {\sf I}(f) \| =\| f\| _{\infty }$, then
$\sf I$ extends to a linear isometry (also called $\sf I$) of 
$L_{\infty }(P)$ onto ${\cal L}(H)$.
\par If $f\in L_{\infty }(P)$, then the operator ${\sf I}(f)$
in ${\cal L}(H)$ is called the spectral integral of $f$ with respect
to $P$ and is denoted 
\par $(vi)$ $\int_X f(x)P(dx):={\sf I}(f)$.
\par Evidently properties $(I-III,V,VI)$ from \S II.11.8 \cite{fell}
are transferable onto the case considered here. 
These and another properties of the spectral integral are as follows.
\par {\bf 5.1. Propositions.} {\it $(I)$. $\int_Xf(x)P(dx)=\int_Xg(x)P(dx)$
if and only if $f$ and $g$ differ only on a $P$-null set.
\par $(II)$. $\int_Xf(x)P(dx)$ is linear in $f$.
\par $(III)$. $\int_Xf(x)g(x)P(dx)=(\int_Xf(x)P(dx))(\int_Xg(x)P(dx))$
for each $f$ and $g\in L_{\infty }(P)$.
\par $(V)$. $\| \int_Xf(x)P(dx)\| =\| f\|_{\infty }$.
\par $(VI)$. If $A\in \sf L$, then $\int_XCh_A(x)P(dx)=P(A)$,
\par in particular $\int_XP(dx)=P(X)={\bf 1}_H$.
\par $(VII).$ For each pair $\xi \in H$ and $\eta ^*\in H^*$, 
let $\mu _{\xi ,\eta }(A):=
\eta ^*(P(A)\xi )$ for each $A\in \sf L$. If $E=\int_Xf(x)P(dx)$
then $\eta ^*(E\xi )=\int_Xf(x)\mu_{\xi ,\eta }(dx)$.
\par $(VIII).$ If $A\in \sf L$, then $P(A)$ commutes with
$\int_Xf(x)P(dx)$, where $e^i:=e_i^*$ such that $e^i(e_j)=\delta ^i_j$.}
\par An $H$-projection-valued measure $P$ on 
the algebra $\sf L$ containing an algebra $Bco(X)$ of clopen (closed 
and open at the same time) subsets of $X$ is called
an $H$-projection-valued tight measure on $X$. We call
$P$ regular if 
\par $(vii)$ $P(A)=\sup \{ P(C): C\subset A\mbox{ and }C \mbox{ is
compact } \} $ for each $A\in \sf L$, where $\sup $ is the least
closed subspace of $H$ containing $range\mbox{ }P(C)$ and to it
corresponds projector on this subspace. Indeed, $P(A)H$ is closed in $H$,
since $P^2(A)=P(A)$. Therefore, 
\par $(viii)$ $P(A)=\inf \{ P(U): U\mbox{ is open and }
U\supset A \} =I- \sup \{ P(C): C\subset X\setminus A\mbox{ and }
C\mbox{ is compact } \} $, hence 
\par $(ix)$ the infimum corresponds to the
projection on $\bigcap_{U\supset A, U\mbox{ is open}}P(U)H$.
A measure $\mu : {\sf L}\to \bf K$ is called regular, if for each 
$\epsilon >0$ and each $A\in \sf L$ with $\| A\|_{\mu }<\infty $ there
exists a compact subset $C\subset A$ such that 
$\| A\setminus C\|_{\mu }<\epsilon $. Since $\| P(X)\| =1$, 
then $\| \mu _{\xi ,\eta }\| \le \| \xi \|_H \| \eta \|_{H^*}$.
For the space $H$ over $\bf K$
measures $\mu _{\xi ,\eta }$ on locally compact $X$ are tight 
for each $\xi , \eta $ in a subset $J\subset H\hookrightarrow H^*$ 
separating points of $H$ if and only if $P$ is defined on 
${\sf L}\supset Bco(X)$; $P$ is regular if and only if
$\mu _{\xi ,\eta }$ are regular for each $\xi , \eta \in J$
due to Conditions $(viii)$ and $(ix)$.
We can restrict our consideration by $\mu _{\xi ,\xi }$ 
instead of $\mu _{\xi ,\eta }$ with $\xi , \eta \in
sp_{\bf K}J$, since ${+\choose -}2\mu _{\xi ,\eta }
=\mu _{\xi {+\choose -}\eta ,\xi {+\choose -}\eta }-\mu _{\xi ,\xi }
-\mu _{\eta ,\eta }$.
\par By the closed support of an $H$-projection-valued tight measure
$P$ on $X$ we mean the closed set $D$ of all those $x\in X$ such that
$P(U)\ne 0$ for each open neighbourhood $x\in U$, $supp\mbox{ }(P):=D$.
\par {\bf 6.} We fix a locally compact totally 
disconnected Hausdorff space $X$
and a Banach space $H$ over $\bf K$ and let $T: C_{\infty }(X,{\bf K})
\to {\cal L}(H)$ be a linear continuous map from the $C$-algebra 
$C_{\infty }(X,{\bf K})$ of functions $f: X\to \bf K$ such that:
\par $(i)$  $T_{fg}=T_fT_g$ for each $f$ and 
$g\in C_{\infty }(X,{\bf K})$,
\par $(ii)$ $T_{\bf 1}=I$ for compact $X$.
\par In general $C_{\infty }(X,{\bf K})$ can be considered as the
$E$-algebra if define $f^t:=f$ for each $f\in C_{\infty }(X,{\bf K})$,
so we can put $T_f^t=T_f$, but the latter equality will not be used.
\par From this definition it follows, that $\| T\| \le 1$, since
$T_{f^n}=T_f^n$  for each $n\in \bf Z$ and $f\in C_{\infty }(X,{\bf K})$.
If $X$ is locally compact and is not compact, 
then $X_{\infty }:=X\cup \{ x_{\infty } \}$
be its one-point Alexandroff compactification. 
Each $f\in C(X_{\infty },{\bf K})$
can be written just in one way in the form $f=\lambda {\bf 1}+g$,
where $g\in C_{\infty }(X,{\bf K})$ and $\bf 1$ is the unit function
on $X_{\infty }$. Therefore, we can extend $T: C_{\infty }(X,{\bf K})
\to {\cal L}(H)$
to a linear map $T': C(X_{\infty },{\bf K})\to {\cal L}(H)$ by setting
${T'}_{\lambda {\bf 1}+g}=\lambda {\bf 1}_H+T_g$ such that 
${T'}_{\bf 1}={\bf 1}_H$.
\par Therefore, $f\mapsto \eta ^*(T_f\xi )=:\tilde \mu _{\xi ,\eta }(f)$ 
is a continuous $\bf K$-linear
functional on $C_{\infty }(X,{\bf K})$, where $\xi \in H$ and $\eta ^*
\in H^*$.
In view of Theorems 7.18 and 7.22 \cite{roo} about correspondence
between measures and continuous linear functionals 
(the non-Archimedean analog of the F. Riesz 
representation theorem) there exists the unique measure
$\mu _{\xi ,\eta }\in M(X)$ such that 
\par $(I)$ $\eta ^*(T_f\xi )=\int_X f(x)\mu _{\xi ,\eta }(dx)$ for each
$f\in C_{\infty }(X,{\bf K})$. 
In the case $T_f^t=T_f$ we have $\mu _{\xi ,\eta }=
\mu _{\eta ,\xi }$, when $\xi ,\eta \in H$. Since
$T_{\bf 1}=I$, then $\mu _{\xi ,\eta }(X)=\eta ^*(\xi )=\xi ^*(\eta )$.
Then for each $A\in \sf L$ we have $\| A\|_{\mu _{\xi ,\eta }}
\le \| \xi \| \mbox{ } \| \eta \| \mbox{ }\sup_{f\ne 0} \| T_f\|$
$\le \| \xi \| \mbox{ } \| \eta \| $. Since $H$ considered as a subspace
of $H^*$ separates points in 
$H$, then for each $A\in \sf L$ there exists the unique linear operator
$P(A)\in {\cal L}(H)$ such that:
\par $(II)$  $\| P(A)\| \le 1$ and  $\eta ^*(P(A)\xi )=\mu _{\xi ,\eta }(A)$,
since $\mu _{\xi ,\eta }(A)$ is a continuous bilinear $\bf K$-valued
functional by $\xi $ and $\eta \in H$.
From the existence of the $H$-projection-valued 
measure in the case of locally compact $X$
we get a projection-valued measure $P'$ on $X_{\infty }$ such that
\par $(III)\quad {T'}_f=\int_{X_{\infty }}f(x)P'(dx)$ 
for each $f\in C(X_{\infty },{\bf K})$.
\par {\bf 7.} {\bf Lemma.} {\it For each $A$ and $B\in \sf L$:
\par $(i)$ $P(A\cap B)=P(A)P(B)=P(B)P(A)$.}
\par {\bf Proof.} For each $g\in C_{\infty }(X,{\bf K})$ 
and $\xi ,\eta \in H$ let
$\nu _g(dx):=g(x)\mu _{\xi ,\eta }(dx)$. For each $f$ and 
$g\in C_{\infty }(X,{\bf K})$ we have: $\int_Xf(x)\mu _{T_g\xi ,\eta }(dx)=
\eta ^*(T_fT_g\xi )=\int_Xf(x)g(x)\mu _{\xi ,\eta }(dx)=\int_Xf(x)\nu _g(dx)$,
consequently, $\nu _g=\mu _{T_g\xi ,\eta }$. For a fixed $A\in \sf L$ let
$\rho (B):=\mu _{\xi ,\eta }(A\cap B)$ for each $B\in \sf L$. Therefore,
$\rho $ is a tight measure on $X$: $\rho \in M(X)$, 
where $M(X)$ denotes the set
of all tight measures on $X$. For each
$g\in C_{\infty }(X,{\bf K})$
there are equalities: 
$$\mbox{ }\int_Xg(x)\rho (dx)=\int_Ag(x)\mu _{\xi ,
\eta }(dx)=\nu _g(A)=\int_Xg(x)\mu _{P(A)\xi ,\eta }(dx).$$ 
Then for each $B\in \sf L$ we get: 
$$\eta ^*(P(A\cap B)\xi )=
\mu _{\xi ,\eta }(A\cap B)=\rho (B)=\mu _{P(A)\xi, \eta }(B)=
\eta ^*(P(A)P(B)\xi ).$$ 
The elements $\xi $ and $\eta \in H$ were arbitrary, hence
$P(A\cap B)=P(A)P(B)$. Interchanging $A$ and $B$ we get the conclusion
of this lemma.
\par {\bf 8.} {\bf Corollary} {\it For each $A\in \sf L$ we have
$P^2(A)=P(A)$ and $P(A)$ is a projection operator such that
$P(X)=I$. If $A\cap B=\emptyset $,
$A$ and $B\in \sf L$, then $P(A)P(B)=0$.}
\par {\bf 9.} {\bf Proposition.} {\it If Conditions $5.(i-iii)$
are satisfied, then $P$ is the unique regular $H$-projection-valued
tight measure on $X$ and $T_f=\int_Xf(x)P(dx)$ for each 
$f\in C_{\infty }(X,{\bf K})$.}
\par {\bf Note.} Such integral is called the spectral integral.
\par {\bf Proof.} Let $\{ A_n: n\in {\bf N} \} $ be a sequence of pairwise
disjoint subsets of $X$, $A_n\in \sf L$. 
Since $X$ is locally compact, then the spectral integral defined in \S 5 
as the limit of certain finite sums exists.
By Corollary 8 $P(A_n)$ are pairwise
orthogonal projectors. Put $Q=\sum_nP(A_n)$. Then for each
$\xi $ and $\eta \in H$ we have: $\eta ^*(Q\xi )=\sum_n\eta ^*(P(A_n)\xi )$
$=\sum_n\mu _{\xi ,\eta }(A_n)=\mu _{\xi ,\eta }(\bigcup_nA_n)$
$=\eta ^*(P(\bigcup_nA_n)\xi )$, consequently,
$P(\bigcup_nA_n)=\sum_nP(A_n)$ and $P$ is an $H$-projection-valued measure.
Since $X$ is locally compact, then each 
measure $\mu _{\xi ,\eta }$ is tight and regular
(see Theorem 7.6 in \cite{roo}), hence $P$ is regular (see \S 5).
Take $f\in C_{\infty }(X,{\bf K})$ and form the spectral integral
$E=\int_Xf(x)P(dx)$. For each $\xi , \eta  \in H$ we have 
$\eta ^*(E\xi )=\int_Xf(x)\mu _{\xi ,\eta }(dx)=\eta ^*(T_f\xi )$,
consequently, $E=T_f$.
In view of Equality $6.(III)$ we have a regular $H$-projection valued measures
both in the case of compact and non-compact locally compact $X$.
\par It remains to verify the uniqueness of $P$. Suppose there exists
another regular $H$-projection-valued tight 
measure on $X$ with the same properties. 
Put $\mu _{\xi ,\eta }(A)=\eta ^*(P(A)\xi )$, $\nu _{\xi ,\eta }
(A)=\eta ^*(Q(A)\xi )$ for each $A\in \sf L$, where $\xi $ and $\eta \in H$.
Then $\int_Xf(x)\mu _{\xi ,\eta }(dx)=\eta ^*(T_f\xi )=
\int_Xf(x)\nu _{\xi ,\eta }(dx)$ for each $f\in C_{\infty }(X,{\bf K})$, 
hence
$\mu _{\xi ,\eta }=\nu _{\xi ,\eta }$ for each
$\xi , \eta \in H$, consequently, $P(A)=Q(A)$ for each $A\in \sf L$.
\par {\bf 10.} {\bf Corollary.} {\it The relation $T_f=\int_Xf(x)P(dx)$ for 
each $f\in C_{\infty }(X,{\bf K})$ sets a one-to-one correspondence between
the set of all regular $H$-projection-valued tight measures
$P$ on $X$ and the set of all continuous linear maps $T:
C_{\infty }(X,{\bf K})\to {\cal L}(H)$, which satisfy conditions 
$5.(i-iii)$.}
\par {\bf 11. } {\bf Note.} A particular case of $H=C_{\infty }(X,{\bf K})$
for locally compact totally disconnected Hausdorff space $X$
and $T_f=f$ for each $f\in C_{\infty }(X,{\bf K})$ 
can be considered independently 
from the given above and it is the following.
Each such $f$ is a limit of a certain sequence by $n\in \bf N$
of finite sums
$\sum_jf(x_{j,n})Ch_{V_{j,n}}(x)$, where $\{ V_{j,n}:
j\in \Lambda _n \} $ is a finite partition of $X$ into the disjoint union
of $V_{j,n}$ clopen in $X$, $x_{j,n}\in V_{j,n}$, 
$\Lambda _n\subset {\bf N}$, since $Range\mbox{ }(f)$ is bounded.
If take $P(V)=Ch_V$ for each $V\in \sf L$, then
$T_ag=\lim_{n\to \infty }\sum_jf(x_{j,n})Ch_{V_{j,n}}(x)g
=\int_Xf(x)P(dx)g$
for each $g\in H$, so there is the bijective correspondence between
elements $a\in \sf A$ of a $C$-algebra $\sf A$ realised as 
$C_{\infty }(X,{\bf K})$
with $X=Sp ({\sf A})$ and their spectral integral representations.
It can be lightly seen that $P(V_1\cap V_2)=Ch_{V_1\cap V_2}=Ch_{V_1}Ch_{V_2}$
$=P(V_1)P(V_2)=P(V_2)P(V_1)$ for each $V_j\in \sf L$. If $\{ V_j:
V_j\in {\sf L}, j\in {\bf N} \} $ is a disjoint family, then $P(\bigcup_jV_j)g=Ch_{\bigcup_jV_j}g=
\sum_jCh_{V_j}g=$ $\sum_jP(V_j)g$ for each $g\in H$. Also $P(\emptyset )H=
Ch_{\emptyset }H=\{ 0 \} $ and $P(X)g=Ch_Xg=g$ for each $g\in H$.
Therefore, $P$ is indeed an $H$-projection-valued tight measure. 
\par Suppose now that $X$ is not locally compact, for example, $X=c_0(\omega _0,
{\bf S})$ with an infinite residue class field $k$ of 
a non-Archimedean infinite field $\bf S$ with non-trivial valuation. 
Then there are $f\in C_{\infty }(X,{\bf K})$ 
for which convergence of finite or even countable or of the cardinality
$card\mbox{ }(k)$ (which may be greater or equal to 
$card \mbox{ }({\bf R})$) sums $\sum_jf(x_{j,n})Ch_{V_{j,n}}$
becomes a problem for a disjoint family $\{ V_{j,n}: j \} $ of 
clopen in $X$ subsets, since $\| Ch_{V_{j,n}} \|_{C(X,{\bf K})}=1$
for each $j$ and $n$.
\par {\bf 12. Theorem. (The non-Archimedean analog of the Stone theorem.)}
{\it  Let $\sf A$ be a commutative Banach $C$-algebra
over a locally compact field $\bf K$. If $P$  is a regular
$H$-projection-valued tight measure on $\hat A$ 
(see \S \S 6,9), then the equation
\par $(i)$ $T_a=\int_{\hat A}\hat a(\phi )P(d\phi )$
for each $a\in \sf A$ defines a 
non-degenerate representation of $\sf A$ in $H$.
Conversely, each non-degenerate representation $T$ of $\sf A$ on a 
Banach space $H$ determines a unique regular $H$-projection-valued tight 
measure $P$ on $\hat A$ such that $(i)$ holds.}
\par {\bf Proof.} The right side of $(i)$ is the spectral integral.
Let $P$ be a regular $H$-projection-valued tight measure on $\hat A$. By
Corollary 10 $T': f\to \int_Xf(x)P(dx)$ is a non-degenerate representation of
$C_{\infty }(\hat A,{\bf K})$ on $H$. By Proposition 5.2 the map 
$a\mapsto \hat a
|_{\hat A}$ is a homomorphism of $\sf A$ onto a dense subset of a subalgebra
of $C_{\infty }(\hat A,{\bf K})$ such that the map 
$T: a\mapsto T'|_{(\hat a|_{\hat A})}=
\int_{\hat A}\hat aP(d\hat a)$ is a non-degenerate representation of $\sf A$.
\par Conversely, let $T$ be a non-degenerate representation of $\sf A$
on $H$. Then from \S 1,2 and Proposition 5.2 it follows, that there
exists a non-degenerate representation $T'$ of 
$C_{\infty }(\hat A,{\bf K})$ such that 
\par $(ii)$ $T_a={T'}_{(\hat a|_{\hat A})}$ for each
$a\in \sf A$. In view of Proposition 9 there exists 
a regular $H$-projection-valued tight measure $P$ on $\hat A$ satisfying
\par $(iii)$ ${T'}_f=\int_{\hat A}f(x)P(dx)$ for each 
$f\in C_{\infty }(\hat A,{\bf K})$.
Combining $(ii,iii)$ we get Formula $(i)$. Let $Q$ be another
regular $H$-projection-valued tight measure which is also related
with the representation $T$ by Formula $(i)$, then 
\par $(iv)$ $\int_{\hat A}
\hat a(x)Q(dx)=T_a=\int_{\hat A}\hat a(x)P(dx)$ for each $a\in \sf A$.
Due to Proposition 5.2 $\{ \hat a|_{\hat A}: a\in {\sf A} \} $
is dense in $C_{\infty }(\hat A,{\bf K})$ 
with respect to the supremum-norm. Hence from 
$(iv)$ and \S 5 it follows that $\int_{\hat A}f(x)P(dx)=\int_{\hat A}
f(x)Q(dx)$ for each $f\in C_{\infty }(\hat A,{\bf K})$, 
consequently, by Proposition 9
$Q=P$.
\par {\bf 13.} {\bf Definition.} $P$ from Theorem 12 is called
the spectral measure of the non-degenerate representation $T$ of $\sf A$.  
\par {\bf 14. Proposition.} {\it Let $P$ be the spectral measure
of the non-degenerate representation $T$ of a commutative Banach 
$C$-algebra $A$ over a locally compact field $\bf K$. 
If $\Omega \subset \hat A$ and $\Omega \in \sf L$, then 
\par $(i)$ $range\mbox{ }(P(\Omega ))=
\bigcup_{\phi \in \Omega }
\{ \xi \in H(T): T_a\xi =\phi (a)\xi \mbox{ for each }a\in {\sf A} \} $.}
\par {\bf Proof.} Relation $12.(i)$ and the definition of the spectral integral
show, that if $\xi \in range\mbox{ }P(V)$ for each $V\in \sf L$ with
$\phi \in V$, then $T_a\xi =\phi (a)\xi $ for each $a\in \sf A$.
\par Conversely, suppose that $T_a\xi =\phi (a)\xi $ for each $a\in \sf A$.
If $T'$ is the representation of $\sf A_t$ isomorphic with $C_{\infty }
(\hat A,{\bf K})$
and $T'$ corresponds to $T$, then 
\par $(ii)$ ${T'}_f\xi =f(\phi )\xi $
for each $f\in C_{\infty }(\hat A,{\bf K})$.\\
Assume $\xi \notin range\mbox{ }(P(\Omega ))$ and consider a measure
$\mu _{\xi ,\eta }(W):=\eta ^*(P(W)\xi )$ for $\xi $ and $\eta \in H$.
There exists $\eta =\xi \ne 0$ such that $\mu _{\xi ,\xi }$ is not
carried by $\Omega $. Due to regularity of $\mu _{\xi ,\xi }$
there exists a compact $E\subset \hat A$, $E\in \sf L$, 
$E\subset \Omega $ such that $\phi \notin
E$ and $\| E \| _{\mu _{\xi ,\xi }}>0$. We take 
$f\in C_{\infty }(\hat A,{\bf K})$
which is not equal to zero everywhere on $E$ and $f(\phi )=0$, since
$\hat A$ is the completely regular topological space $T_{3.5}$ 
(see Theorem 2.3.11 in \cite{eng}).
From Formula $(ii)$  it follows, that ${T'}_f\xi =0$.
By Chapter 7 in \cite{roo} and Formula $5.1.(VII)$ above
there is an inequality:
$\| {T'}_f\xi \| \ge \| f \| _{N_{\mu _{\xi ,\xi }}}\ge \sup_{x\in E}
|f(x)|N_{ \mu _{\xi ,\xi }}(x)=:\| f|_E \| _{N_{\mu _{\xi ,\xi }}}$,
where $\| f\|_{\phi }:=\sup_{x\in X} |f(x)|\phi (x)$ for $f: X\to \bf K$ 
and $\phi : X\to [0,\infty )$;
$N_{\mu }(x):=\inf_{U\in {\sf L}, x \in U}\| U\| _{\mu }$;
$\| A\|_{\mu }:=\sup \{ |\mu (B)|:
B\in {\sf L}, B\subset A \} $ for each $A\in \sf L$.
If $\| \xi \|_H=1$, then $\| {T'}_f\xi \|=\| f\|_{N_{\mu _{\xi ,\xi }}}$.
We get a contradiction, consequently, $\xi \in range\mbox{ }P(\Omega )$.
\par {\bf 15.} Let $\sf A$
be a commutative $C$-algebra with the unit $\bf 1$ 
over a locally compact field $\bf K$ and let $\mu $
be any regular tight measure on $\hat A$. 
Let the space $L(\hat A,\mu ,{\bf K})$ be defined on the algebra
$\sf L$ such that ${\sf L}\supset Bco(\hat A)$ 
of $\hat A$ as in Chapter 7 \cite{roo}:
it is the completion relative to $\| f\|_{N_{\mu }}$
of the $\bf K$-linear space of all step functions, that is,
finite linear combinations of characteristic functions
of elements of $\sf L$.
\par {\bf 15.1. Theorem.} {\it The equation 
\par $(i)$ $(T_af)(\phi )=\hat a
(\phi )f(\phi )$ for each $a\in \sf A$, $f\in L(\hat A,\mu ,{\bf K})$
and $\phi \in \hat A$ defines a non-degenerate representation $T$ of $\sf A$
on $H=L(\hat A,\mu ,{\bf K})$ and the spectral measure $P$ of $T$
is given by $P(W)f=Ch_Wf$ for each $W\in \sf L$ and $f\in H$.}
\par {\bf Proof.} If $\hat a(\phi )\in C(\hat A)$, 
then $\sup_{\phi \in \hat A}
|\hat a(\phi )|<\infty $, so 
$$\sup_{\phi \in \hat A}|f(\phi )|
\mbox{ }|\hat a(\phi )|\mbox{ }N_{\mu }(\phi )\le \| f\|_{\mu }
\| \hat a\|_{C(\hat A,{\bf K})},$$ 
consequently, $\hat af\in H$
and $T_af=\int_{\hat A}\hat a(\phi )P(d\phi )f$  due to Theorem 12.
Therefore, for $\hat a=Ch_W$ we have $T_af=P(W)f$ for each $W\in Bco(\hat A)$.
Each measure
$\mu _{\xi ,\eta }(W)=\eta ^*(P(W)\xi )$ has an extension from $Bco(\hat A)$
on $\sf L$ due to its regularity, 
where $\xi $ and $\eta \in H$ (see \S 5). 
The family of such measures $\mu _{\xi ,\eta }$ characterise $P$ completely, 
since $H$ is the Banach space of separable type over $\bf K$.
Therefore, we have an extension of $P$ on $\sf L$. 
\par {\bf 16.} From the results above
it follows that $supp\mbox{ }(P)\subset Sp({\sf A})$.
If a representation $T$ is one to one, that is, $ker\mbox{ }T=\{ 0\} $,
then $T$ is called faithful. 
\par {\bf 16.1. Proposition.} 
{\it The kernel of a non-degenerate representation
$T$ of $\sf A$ consists of $a\in \sf A$ such that $\hat a$ vanishes everywhere
on the spectrum of $T$. Suppose in addition that $\sf A$ is a commutative
$C$-algebra over a locally compact field $\bf K$. 
Then $T$ is faithful if and only if its spectrum is all
of $\hat A$.}
\par {\bf Proof.} A condition $T_a=0$ is equivalent to 
$\int_{\hat A}\hat a(\phi )
P(d\phi )=0$, which is equivalent to $\hat a(\phi )|_{\hat A}=0$
due to Theorem 12.
Therefore, $ker \mbox{ }T=\{ 0\} $
is equivalent to $supp\mbox{ }T=\hat A$.
\par {\bf 17.} Fix a Banach space $H$ over a non-Archimedean complete 
field $\bf F$ such that ${\bf F}\subset \bf C_p$.
If $b\in {\cal L}(H)$ we write shortly $Sp(b)$ instead of
$Sp_{{\cal L}(H)}(b):=cl(Sp(sp_{\bf F}\{ b^n: n=1,2,3,... \} )) $ 
(see also \cite{roo}). 
If $A$ is a commutative Banach subalgebra in ${\cal L}(H)$,
then there exists a quotient mapping $\theta : A\to A/B_A$, where $B_A$ 
is a closed subalgebra of $A$ such that $B_A=ker (\| * \|_{sp})$
is the kernel of the spectral norm, $B_A:=\{ x: x \in A; 
\| x \| _{sp }=0 \} $. Then $\theta (A)$ is the normed algebra, 
$\theta (A)$ is the subalgebra of ${\cal L}(H)/B_{{\cal L}(H)}$.
Choose a locally compact subfield $\bf K$ in $\bf F$.
\par {\bf 17.1. Spectral theorem for operators.} 
{\it Let $b\in {\cal L}(H)$. 
Then there exists a unique
$H$-projection-valued tight measure $P$ on $\bf K$ with values in
${\cal L}(H)$ with the following properties:
\par $(i)$ the closed support $D$ of $P$ is bounded in $\bf K$;
\par $(ii)$ $\theta (b)=\int_{\bf K}xP(dx);$
also $b=\lambda _b V\int_{\bf K}xP(dx)$, where 
$V$ is a continuous operator from $H$ onto its closed 
$\bf K$-linear subspace such that $|\pi _{\bf K}|
\le \| V\| \le |\pi _{\bf K} | ^{-1}$, $\pi _{\bf K}\in \bf K$,
$|\pi _{\bf K}|=\max \{ |x|: x\in {\bf K}, |x|<1 \} $;
$|\lambda _b|=\| b \| $, $\lambda _b\in \bf F$;
$V$ is an isometry of $H$ onto its closed $\bf K$-linear subspace
for ${\bf K}=\bf F$;
\par $(iii)$ if ${\bf K}=\bf F$, then $D=Sp(a)$, where
$a$ is an auxiliary operator defined by $V$ and $b$.
\par Moreover, if $S$ is a family of commuting operators,
$S\subset {\cal L}(H)$, then there exists a unique $H$-projection-valued 
tight measure $P$ on a locally compact subset $X\subset 
B({\bf K},0,1,)^{\gamma }$ such that for each $b\in S$ 
there exists a unique $f_b\in C_{\infty }(X,{\bf F})$ for which
$\theta (b)=\int_Xf_b(x)P(dx)$, 
$b=V \int_Xf_b(x)P(dx)$ and $V$ as above,
where $cl\mbox{ }X=X\cup \{ 0 \} $ is compact.}
\par {\bf Proof.} If $\| b \|_{sp}=0$, then $Sp (b)=\emptyset $ and this 
case is trivial with $P(\emptyset )=0$. 
In the case of 
the locally compact field $\bf K$ and $H$ over $\bf K$ 
we can take $W:=cl(b(H))$ and $b=Va$, where $V$ is an operator
such that $V(H)=W$, $V|_W: W\to W$ is an isometry,
$V(H\ominus W)=\{ 0 \} $, $a$ is an operator of $H$ onto $H$
such that $\| a \|_{sp}>0$ and $a$ is representable 
as a convergent series of projectors in some basis
of $H$, that can be shown by transfinite induction
(see \cite{put}). In this case we get $(ii)$.
Analogosuly for commuting algebras of operators.
\par The field $\bf F$ can be considered as the Banach 
space over $\bf K$. This means that $\bf F$ supplied with the linear 
structure over $\bf K$ is isomorphic with $c_0(\beta ,{\bf K})$
and a corresponding ordinal $\beta $, since
$\bf K$ is locally compact and hence spherically complete (see Theorems
5.13 and 5.16 \cite{roo}). This isomorphism $\chi :
{\bf F}\to c_0(\beta ,{\bf K})$ may be non-isometrical.
The isomorphism $\chi $ generates the isomorphism of $H$ considered
as the Banach space $H_{\bf K}$ over $\bf K$ with 
$c_0(\alpha _{\bf K},{\bf K})$, $\chi :H_{\bf K}\to 
c_0(\alpha _{\bf K},{\bf K})$ with the corresponding ordinal 
$\alpha _{\bf K}$. This isomorphism is $\bf K$-linear
and it produces an injective continuous $\bf K$-linear
embedding $\chi ^*: {\cal L}(H)\to {\cal L}
(c_0(\alpha _{\bf K},{\bf K}))$ with continuous
$(\chi ^*)^{-1}|_{\chi ^*({\cal L}(H))}$. 
The embedding $\chi ^*$ is given by the following formula:
$\chi ^*(a)y:= \chi a \chi ^{-1}y$ for each $a\in {\cal L}(H)$ 
and $y\in c_0(\alpha _{\bf K},{\bf K})$. This is the 
well-known construction of the contraction of a scalar field
for a Banach space. In the particular case of
$H=\bf F^n$ with $n\in \bf N$ and if $\bf F$ is a finite algebraic
extension of $\bf K$ to each $a\in {\cal L}(H)$ there corresponds a 
finite $n\times n$ matrix, hence 
$\chi ^*(a)\in {\cal L}({\bf K}^{\beta n})$.
\par Suppose there is a representation of a $C$-algebra
$C_{\infty }(X,{\bf F})$ with the help of a 
$c_0(\alpha _{\bf K},{\bf K})$-projection-valued tight measure $P$ on
a locally compact subset $X$ in ${\bf K}^{\gamma }$.
Then $(\chi ^*)^{-1}$ produces from $P$ an 
$H$-projection-valued tight measure $P_{\bf F}$, 
since
$$(\chi ^{-1}P(V)\chi )(\chi ^{-1}P(W)\chi )=
\chi ^{-1}P(V\cap W)\chi $$ 
for each $V$ and $W\in \sf L$. Consequently, 
$$\chi ^{-1}\int_Xg(x)P(dx)\chi z=\int_X(\chi ^{-1}g(x)
\chi )(\chi ^{-1}P(dx)\chi )z$$ 
for each $z\in H$ and 
$f\in C_{\infty }(X,{\bf K})$.
Therefore, if $\chi ^*(a)=\int_Xg_a(x)P(dx)$, then
$a=\int_Xf_a(x)P_{\bf F}(dx)$, where $g_a\in C_{\infty }(X,{\bf K})$ and
hence $\chi ^*(g_a)=:f_a\in C_{\infty }(X,{\bf F})$
such that $f_a=g_a$, since the restriction of $\chi $ on $\bf K$
embedded into $\bf F$ is the identity $\bf K$-linear mapping.
\par From this it follows, that instead of $b$ or $S$ it is
sufficient to consider $\chi ^*(b)$ or $\chi ^*(S)$.
Denote $\chi ^*(b)$ and $\chi (S)$ simply by $b$ and $S$ respectively.
The operator $b$ or the family $S$ generates a commutative subalgebra
$A$ of ${\cal L}(H_{\bf K})$ generated by 
$sp_{\bf K}\{ b^n: 0\le n\in {\bf Z} \} $
or by $sp_{\bf K}\{ a_1^{m_1}...a_n^{m_n}:$ $a_j\in S, n\in {\bf N},
0\le m_j\in \bf Z;$ $j=1,...,n \} $, where $b^0:=\bf 1$.
It has a completion $\sf A$ relative to the spectral norm $\| *\|_{sp}$
by Chapter 6 in \cite{roo}.
This $\sf A$ in view of Theorem 6.15 and Corollaries 6.16, 6.17
in \cite{roo}
is a Banach $C$-algebra $C(Sp({\sf A}),{\bf K})$, 
since ${\bf 1}\in {\sf A}$ and $Sp(A)$ is compact such that
$C(E,{\bf K})$ is isomorphic with $C_{\infty }(E,{\bf K})$
for compact $E$.
\par For $S=\{ b \} $ each $\phi \in A^+$ is completely defined by
$\phi (b)$, thus the map $\phi \mapsto \phi (b)$ is continuos and one-to-one, consequently, it is a homeomorphism from the compact space $Sp({\sf A})$
onto the compact subset $Sp(b)$ of $\bf K$. Therefore, we identify
$Sp({\sf A})$ with $Sp(b)$. 
So the Gelfand transform of $b$ becomes the identity function on
$Sp(b)$. Thus by Theorem 12 the identity representation of $\sf A$ in
$H$ gives rise to the $H$-projection-valued tight measure $P$ on $Sp(b)$
such that $\theta (b)=\int_{D}xP(dx)$
and $b=\lambda _bV \int_{D}xP(dx)$. 
Since the identity representation is faithful, 
then Proposition 16 shows that the closed support $D$ of $P$ is homeomorphic
to $Sp(b)$ (up to the mapping $\chi |_{\bf K}$). 
Extending $P$ to $\sf L$ by setting $P(W)=P((W)\cap
Sp(b))$, we thus obtain a tight projection-valued measure on $\bf K$ 
satisfying $5.(i-iii).$ 
\par In the case of the family $S$ take $X=Sp({\sf A})$
such that due to Theorem 2.3.20 in \cite{eng} about a diagonal mapping
we can choose ${\sf a}\le card\mbox{ }(\gamma )\le 
card ({\sf A}^+)$ while embedding $X\hookrightarrow B({\bf K},0,1)$,
where $\sf a$ is the minimal cardinality of a family
of subsets in ${\sf A}^+$ separating points of $\sf A$.
In view of \S 6.2 in \cite{roo} $cl(X)=X\cup \{ 0\} $.
\par To show that $P$ is uniquely determined by
$(i,ii)$, let $P'$ be another $H$-projection-valued tight measure
on $\bf K$ satisfying $(i,ii)$ and let $E$ be a compact subset of $\bf K$
containing the supports of both $P$ and $P'$. Consider two representations
$T: f\mapsto \int_Ef(x)P(dx)$ and $T': f\mapsto \int_Ef(x)P'(dx)$
of $C(E,{\bf K})$. If $w$ and $e$ are the identity function $w(x)=x$ and 
the constant function $e(x)=1$ for each $x\in E$ respectively, then
condition $(ii)$ satisfied by both $P$ and $P'$ shows that $T_w={T'}_w=b$
and also $T_e={T'}_e=\bf 1$. By the Kaplansky theorem $w$ and $e$ generate
$C(E,{\bf K})$ as a $C$-algebra \cite{roo,sch1}, hence $T'=T$ and by the uniqueness statement of Theorem 12 we have $P'=P$.
\par {\bf 18.} The $P$ of the above theorem is called the spectral measure
of the operator $b$ or of a family $S$. In particular $S$ may be
a commutative subalgebra of ${\cal L}(H)$.
Evidently, each nilpotent operator $v$ in ${\cal L}(H)$ has 
$\| v \|_{sp}=0$. If $v\in {\cal L}_0(H)$ has $e^*_j(ue_i)=0$
for each $j\le i$, then $v$ is nilpotent.
Therefore, $\theta ({\cal L}_0(H))$ is isomorphic with the algebra of 
diagonal operators on $H$. Its spectrum was found in \S 4.
\par {\bf 19.} {\bf Note.} It is an interesting 
property of the non-Archimedean case that a 
condition of a normality of an operator $b$ is not necessary
in Theorem 17.1 apart from the classical case. 
It is not so suprising if remind, that orthogonality
in the non-Archimedean case is quite another than in the classical case.
For example, two vectors $(1,0)$ and $(1,1)$ are orthonormal in $\bf K^2$,
but are not orthogonal in $\bf C^2$. 
\par This work was started by S. Ludkovsky. All matters of this paper 
were thoroughly discussed and investiagted by both authors.
Sections 2-4 were written by B. Diarra and section 5 was written by
S. Ludkovsky.


\begin{thebibliography}{99}
\bibitem{ack} S. Albeverio, R. Cianci and  A. Khrennikov. "On the spectrum of
the $p$-adic position operator".
J. Phys. A : Math. Gen. {\bf 30} (1997), 881-889.
\bibitem{abpk} S. Albeverio, J.M. Bayod, C. Perez-Garcia, R. Cianci and 
A. Khrennikov. "Non-archimedian
analogues of orthogonal and symmetric operators and $p$-adic
quantization". Acta Applicandae Mathematicae. {\bf 57} (1999), 1-33.
\bibitem{ber} V.G. Berkovich. "Spectral theory and analytic geometry 
over non-Archimedean fields". Math. Surveys and Monographs of A.M.S., 
N {\bf 33} (1990).
\bibitem{diar3} B. Diarra. "Ultraproduits ultrametriques de
corps values". Ann. Sci. Univ. Clermont II, S\'er. Math., {\bf 22}
(1984), 1-37. 
\bibitem{diar1} B. Diarra. "Remarques sur les matrices orthogonales
(resp. sym\'etriques) \`a coefficient $p$-adiques".
Ann. Sci. Univ. "Blaise Pascal", Clermont II, S\'er. Math., Facs.
{\bf 26} (1990), 31-50.
\bibitem{diar2} B. Diarra. "An operator on some ultrametric Hilbert spaces".
J. Analysis  {\bf 6} (1998), 55-74.
\bibitem{dun} N. Dunford, J.T. Schwartz. "Linear operators"
(New York, Intersci. Publ., v. {\bf 1} 1958, v. {\bf 2} 1963).
\bibitem{eng} R. Engelking. "General topology" (Moscow, Mir, 1986).
\bibitem{esc1} A. Escassut. "The ultrametric spectral theory".
Periodica Mathem. Hung. {\bf 11} (1980), 7-60.
\bibitem{esc2} A. Escassut, N. Manetti. "Spectral norm of a $p$-adic
Banach algebra". Bull. Belg. Math. Soc. "Simon Stevin"
{\bf 5} (1998), 79-91.
\bibitem{esc3} A. Escassut. "Analytic elements in $p$-adic analysis"
(Singapore, World Scientific, 1985).
\bibitem{fell} J.M.G. Fell,R.S. Doran. "Representations of $*$-algebras,
localy compact groups, and Banach $*$-algebraic bundles" 
(Boston, Acad. Publ., 1988).
\bibitem{gant} F.R. Gantmaher. "The theory of matrices"
(Moscow, Nauka, 1988; New York, Chelsea Publ. Comp., 1977).
\bibitem{gogro} M. Goto, F.D. Grosshans. "Semisimple Lie algebras"
(New York, Marcel Dekker, Inc., 1978).
\bibitem{grus} L. Gruson. "Th\'eorie de Fredholm $p$-adique".
Bull. Soc. Math. France. {\bf 94} (1966), 67-95.
\bibitem{hew} E. Hewitt, K.A. Ross. "Abstract harmonic analysis"
(Berlin, Springer, 1979).
\bibitem{ko1} H.A. Keller, H. Ochsenius A. 
"Algebras of bounded operators on
non-classical orthomodular spaces". Inter. Jour. Theor. Phys. {\bf 33}
(1994), 1-11.
\bibitem{ko2} H.A. Keller, H. Ochsenius A. "Bounded operators on
non-Archimedean spaces". Math. Slovaca. {\bf 45} (1995), 413-434.
\bibitem{khrb} A. Khrennikov. "Non-Archimedean analysis: Quantum paradoxes,
Dynamical systems and Biological models" (Dordrecht, 
Kluwer Acad. Publ., 1997).
\bibitem{nai} M.A. Naimark. "Normed rings" (Moscow, Nauka, 1968). 
\bibitem{put} M. van der Put. "The ring of bounded operators
on a non-Archimedean normed linear space".
Indag. Math. {\bf 71: 3} (1968), 260-264.
\bibitem{put2} M. van der Put. "Difference equations over $p$-adic fields".
Math. Ann. {\bf 198} (1972), 189-203. 
\bibitem{reed} M. Reed, B. Simon. "Methods of modern mathematical physics"
(New York, Acad. Press, 1975).
\bibitem{ocsc} H. Ochsenius, W.H. Schikhof. "Banach spaces over 
fields with on infinite rank valuation". In: "$p$-Adic Functional 
Analysis". Edited by J. K\c akol, N. De Grande-De Kimpe , C. Perez-Garcia
(New-York, M. Dekker Inc. 1999), 233-293.
\bibitem{roo} A.C.M. van Rooij. "Non-Archimedean functional analysis"
(New York, Marcel Dekker Inc., 1978).
\bibitem{schroo} A.C.M. van Rooij, W.H. Schikhof.
"Non-Archimedean commutative $C^*$-algebras".
Indag. Math. {\bf 35} (1973), 381-389.
\bibitem{sch1} W.H. Schikhof. "Ultrametric calculus" (Cambridge,
Cambr. Univ. Press, 1984).  
\bibitem{vish} M.M. Vishik. "Non-Archimedean spectral theory".
J. Sov. Math. {\bf 30} (1985), 2513-2554.
\end{thebibliography}
\end{document}